\documentclass[reqno,11pt]{amsart}
\usepackage[dvipsnames]{xcolor}
\usepackage{amsmath}
\usepackage{amsfonts}
\usepackage{amssymb}
\usepackage{amsthm}
\usepackage{amsxtra}
\usepackage{mathtools}	
\mathtoolsset{centercolon}		
\usepackage{empheq}
\numberwithin{equation}{section}

\usepackage{bm}					 	
\usepackage[shortlabels]{enumitem}			

\usepackage{esint}					

\usepackage{graphicx}
\usepackage{tikz, calc}
\usetikzlibrary{hobby}
\usepackage{mathrsfs}

\newtheorem{theorem}{Theorem}[section]

\newtheorem{lemma}[theorem]{Lemma}
\newtheorem{lem}[theorem]{Lemma}

\newtheorem{corollary}[theorem]{Corollary}
\newtheorem{assump}[theorem]{Assumption}

\theoremstyle{definition}

\theoremstyle{remark}
\newtheorem{remark}[theorem]{Remark}

\newcommand{\name}[1]{\textsc{#1}}
\newcommand{\df}[1]{\textit{#1}}

\AtBeginDocument{
	\let\div\relax 
	\DeclareMathOperator{\div}{div}
	\let\L\relax
	\newcommand{\L}{\mathcal{L}}
}

\newcommand{\R}{\mathbb{R}}
\newcommand{\Rz}{\mathbb{R}}

\newcommand{\N}{\mathbb{N}}

\newcommand{\Z}{\mathbb{Z}}
\newcommand{\C}{\mathbb{C}}
\newcommand{\Hz}{\mathbb{H}}

\newcommand{\A}{\mathcal{A}}
\newcommand{\Cc}{C}
\newcommand*\dd{\mathop{}\!\mathrm{d}}

\newcommand{\la}{\langle}
\newcommand{\ra}{\rangle}

\newcommand{\epsi}{\varepsilon}
\newcommand{\eps}{\varepsilon}

\newcommand{\ove}{\overline} 
\newcommand{\ti}{\widetilde}

\newcommand{\D}{\nabla}

\newcommand{\id}{\textrm{id}}

\DeclareMathOperator{\SO}{SO}

\newcommand{\wto}{\rightharpoonup} 		
\newcommand{\acto}{\stackrel{\text{\rm \tiny AC}}{\to}}

\newcommand{\loc}{\mathrm{loc}}
\newcommand{\comp}{\mathrm{c}}

\usepackage{caption}
\usepackage{subcaption}

\makeatletter
\def\p@subfigure{\thefigure\,}
\def\p@subtable{\thetable\,}
\makeatother

\usepackage{todonotes}
\usepackage{hyperref}

\title[Linearization in atomistic dynamics]{ Discrete-to-continuum linearization\\
  in atomistic dynamics} 

\begin{document}
	

	\author[M. Friedrich]{Manuel Friedrich} 
\address[Manuel Friedrich]{Department of Mathematics, Friedrich-Alexander Universit\"at Erlangen-N\"urnberg. Cauerstr.~11,
D-91058 Erlangen, Germany, \& Mathematics M\"{u}nster,  
University of M\"{u}nster, Einsteinstr.~62, D-48149 M\"{u}nster, Germany.
}	
	\email{manuel.friedrich@fau.de}
	\urladdr{https://www.math.fau.de/angewandte-mathematik-1/mitarbeiter/prof-dr-manuel-friedrich/}

	\author[M. Seitz]{Manuel Seitz} 
	\address[Manuel Seitz]{Vienna School of Mathematics, University of
		Vienna, Oskar-Morgenstern-Platz~1, A-1090 Vienna, Austria, \& Faculty of Mathematics,  University of
		Vienna, Oskar-Morgenstern-Platz~1, A-1090 Vienna, Austria.
	}
	\email{manuel.seitz@univie.ac.at}

	\author[U. Stefanelli]{Ulisse Stefanelli} 
	\address[Ulisse Stefanelli]{Faculty of Mathematics, University of
		Vienna, Oskar-Morgenstern-Platz 1, A-1090 Vienna, Austria,
		Vienna Research Platform on Accelerating
		Photoreaction Discovery, University of Vienna, W\"ahringerstra\ss e 17, 1090 Wien, Austria,
		\& Istituto di
		Matematica Applicata e Tecnologie Informatiche {\it E. Magenes}, via
		Ferrata 1, I-27100 Pavia, Italy.
	}
	\email{ulisse.stefanelli@univie.ac.at}
	\urladdr{http://www.mat.univie.ac.at/$\sim$stefanelli}
	
	
	
	\begin{abstract}  In the stationary case,  
 atomistic interaction energies can be proved to $\Gamma$\nobreakdash-converge  to   classical elasticity models in the  simultaneous  atomistic-to-continuum  and  linearization
 limit \cite{solci,Schmidt_atomisticToCont}.  The aim of this
note is that of extending the  convergence  analysis  to the dynamic setting.  Moving within the
 framework of \cite{Schmidt_atomisticToCont}, we prove that  solutions of the equation of motion  driven by   atomistic deformation  energies  converge to the
 solutions of the  momentum  equation for the  corresponding continuum  energy of  linearized elasticity.
  By recasting the evolution problems in their equivalent  energy-dissipation-inertia-principle  form, we directly argue at the variational level of evolutionary $\Gamma$\nobreakdash-convergence \cite{Mielke,ss}. This in particular ensures the pointwise in  time convergence of the energies.  
	\end{abstract}
	

	\maketitle

	\pagestyle{myheadings}

	\section{Introduction} 
	
In recent years, the passage from atomistic models to continuum
theories has 
 attracted remarkable attention, originating a variety of
different results, 
see \cite{blancAtomisticContinuumLimits2007,
  braidesDiscretetoContinuumVariationalMethods} for an overview. 
The interest in atomistic-to-continuum limits is grounded on the
possibility of describing macroscopic behaviors starting from
microscopic modeling. 
The latter relies in specifying interactions between displaced points, resulting
from a deformation of a given reference
lattice. The atomistic-to-continuum limit then ensues  by letting the reference lattice spacing
go to  $0$. 

The elective tool to   tackle  these limits is $\Gamma$\nobreakdash-convergence
\cite{braidesGammaconvergenceBeginners2002,
  dalmasoIntroductionGConvergence1993}, see for instance  
\cite{Alicandro_Cicalese_DiscToCont,
  braidesGelliFromDiscreteSystems2006}. Alternatively, more pointwise
perspectives can be followed \cite{blancMolecularModelsContinuum2002},
which again may contribute   to compute   $\Gamma$-limits
in specific cases. 
 These  tools were successfully applied to  different 
situations, ranging from  the modeling of  
lower-dimensional structures
\cite{braidesHomogenizationDiscreteThin2021, Schmidt_plates,
  schmidtPassageAtomicContinuum2008},  to  analyzing fracture
and softening effects in one dimension
\cite{braidesVariationalFormulationSoftening1999,
  braidesContinuumLimitsDiscrete2002, BLO06,
  schaffnerLennardJonesSystemsFinite2017},  to  deriving
cleavage laws in higher dimensions
\cite{friedrichAtomistictoContinuumAnalysisCrystal2014,
  friedrichAnalysisCrystalCleavage2015b,
  friedrichDiscretetocontinuumConvergenceResult2014},  to 
investigating systems of charged particles
\cite{vanmeursDiscretetocontinuumConvergenceCharged2022} and
multi-body systems \cite{bachDiscretetocontinuumLimitsMultibody2019,
  braunPassageAtomisticSystems2013}. Moreover, surface effects
\cite{BLO06, braidesDiscToConPlanarLatticeEnergies,
  delninDiscreteIsoperimetric,
  theilSurfaceEnergiesTwodimensional2011}, as well as random point
clouds and stochastic lattices
\cite{alicandroIntegralRepresentationResults2011,
  braidesAsymptoticBehaviorDirichlet2022, rufApproxMumfordShah} were
 also  studied using these techniques.

The main reference for our work is the stationary analysis by {\sc Braides, Solci,
  \& Vitali} \cite{solci} and {\sc Schmidt}
\cite{Schmidt_atomisticToCont} who recover elasticity models
as $\Gamma$-limits  of 
atomistic energies. In both papers, the atomistic-to-continuum  limit
is combined with a small-deformation limit, effectively resulting in a
linear elastic  limiting  model. We refer to this combined
limiting procedure as {\it atomistic-to-continuum linearization limit}
in the sequel. 

 Following the  general setting of \cite{Schmidt_atomisticToCont}, the focus of
this note is to extend the atomistic-to-continuum linearization limit procedure to the evolutionary case. We  first show that 
solutions to  the equation of motion  driven by   the gradient of the atomistic deformation energy   converge to  the  solution to
 the momentum equation  for the limiting linear elastic energy by passing to the limit in  the equations.  By resorting to an equivalent variational formulation of the evolutionary  problems in terms of the Energy-Dissipation-Inertia Principle, we present an evolutionary $\Gamma$\nobreakdash-convergence analysis \cite{Braides14, Mielke,ss}, as well. This entails in particular the pointwise in time convergence of the energies and the strong convergence of the gradients.

Let us mention that 
atomistic-to-continuum linearization in the evolutionary setting has
already attracted some attention.  The one-dimensional case is studied by {\sc Ortner}~\cite{Ortner_GradFlowAsSelection} who proves that the
 limits of $H^1$-gradient flows of the
Lennard-Jones potential give rise to $H^1$-gradient flows of an elastic
energy. 

 The one-dimensional model is also considered by {\sc Braides, Defranceschi, \&
  Vitali} \cite{braidesVariationalEvolutionOnedimensional2014}, who check that $L^2$-gradient flows of
the Lennard-Jones potential under a different scaling  give rise to
gradient flows of the Mumford-Shah (or Griffith fracture)
functional.  Their analysis is based on an analogous limiting
procedure at the level of time-discrete minimizing-movements
approximations.  Compared with these contributions, the  novelties here are that the multidimensional case is tackled and inertial effects are taken into account.

 Some instances of discrete-to-continuum  limits in dynamics have already been considered.  
 \name{E}  $\&$  \name{Ming} \cite{eCauchyBornRuleStability2007}  obtain a dynamic  atomistic-to-continuum  convergence   result, in which they show that, as long as the deformation gradient stays inside a stability region, solutions to the discrete dynamical problem obtained by molecular mechanics are close to the nonlinear continuum problem obtained by the Cauchy-Born rule. 
 A related result has been achieved by   \name{Ortner}  $\&$
 \name{Theil} \cite{ortnerJustificationCauchyBorn2013}, who prove
 existence of solutions to the discrete and the nonlinear continuum
 problem for more general interaction energies and unbounded domains,
 and likewise quantify the  difference in terms of the lattice
 spacing. Due to the nonlinear nature of the limiting problems and the
 absence of viscous regularization,  in
 \cite{eCauchyBornRuleStability2007} and
 \cite{ortnerJustificationCauchyBorn2013} existence of solutions can
 be  guaranteed for short times only. Existence of global solutions in the bounded-domain case,
 as well as their asymptotic behavior, has
been investigated by {\sc Braun} \cite{Braun}. In all mentioned
 results
 \cite{Braun,eCauchyBornRuleStability2007,ortnerJustificationCauchyBorn2013},
 regularity or smallness assumptions  have to be imposed on data. On the contrary, our approach is nonperturbative,  thus requiring less restrictions on the data.

Let us now describe our result 
more closely,  postponing however  most details  to Section \ref{sec:model}.  We indicate by  $\Omega  \subset {\mathbb
R}^d$ the macroscopic reference configuration of the
specimen. At the atomistic level, the system consists of all points in $\Omega$ 
of the   scaled $d$-dimensional  Bravais  lattice $\epsi\L :=
\epsi A \cdot  \Z^d$, where 
$\epsi > 0$  models   the  interatomic distance  and $A \in
\mathbb{R}^{d\times d}$ is an invertible matrix. The  atomistic
 energy  $I_\epsi(y)$  of a {\it lattice deformation} $y: \Omega \cap \epsi\L \to
\mathbb{R}^d$  
 is  modeled  by  
\begin{align*}
	 I_\epsi  (y) = \epsi^d \det A \sum_{\textrm{cells}}  W 
  ( \bar{\D} y) + \textrm{surface terms}.
\end{align*}
 The latter features the {\it discrete deformation gradient} $\bar{\D} y$, describing the relative displacement of the
atoms of  a  single cell, see \eqref{eq:def of disc grad} below.  
The  sum ranges over all scaled cells of the lattice portion  $\Omega \cap \epsi\L$.
 The scaling $\epsi^d \det A$ corresponds to the volume of 
a single cell  and the function  $W$  models the {\it elastic cell
  energy}, assumed to be suitably smooth.  
The treatment of surface energy  terms  is in general a delicate matter,
see \cite{Schmidt_atomisticToCont}, and, for the sake of conciseness,
will be simplified in this paper by  choosing
homogeneous Dirichlet boundary conditions.  

 In the small-deformation setting, a second small parameter $\delta>0$ is
involved,  quantifying  the distance of deformations from the
identity. More precisely, the atomistic model is rephrased  in terms of
the {\it scaled lattice displacement} $u := \delta^{-1} (y-\id )  : \Omega \cap \epsi\L
\to \mathbb{R}^d$ \cite{DalMaso_Negri_Percivale} and the energy is
suitably rescaled as  
\begin{align*}
	 I_{ \epsi\delta }(u)  = \frac{\epsi^d \det A}{\delta^2}
  \sum_{\textrm{cells}}  W  (  Z  + \delta   \bar{\D} u  
) + \textrm{surface terms},
\end{align*}
 where $Z$ denotes the discrete deformation gradient of $\id$. 
 The main result in  \cite{Schmidt_atomisticToCont} is the
$\Gamma$\nobreakdash-convergence of  $I_{\epsi\delta}$ as  $\epsi,\,\delta \to 0$  to the classical
linear elastic energy 
\[  I( w)  =   \frac12 \int_\Omega  \D^s w:\C :\D^s w \dd x, \] 
 where  $w\colon\Omega \to \mathbb{R}^d$ is now
the
{\it macroscopic displacement} of the body,  $\D^s w := (\D w + \D w ^\top)/2$ is its symmetrized gradient, and
$\C$ the  fourth-order \emph{elasticity  tensor}, which can be directly computed from  $W$, see \eqref{eq:tensor} below. 

 Our  
main result  concerns the limiting behavior as $\epsi, \, \delta
\to 0$ of the  
solutions to the atomistic  equation of motion  
\begin{subequations}	    \label{eq:intro_gf0} 
	\begin{empheq}[left=\empheqlbrace]{align}
		&  \rho \ddot u_{\epsi\delta}(t) + \nu \dot u_{\epsi\delta} (t)   + \partial I_{\epsi\delta} (u_{\epsi\delta}(t))=0
	          \quad \forall t \in (0,T),   \label{eq:intro_gf0_1}\\
		&u_{\epsi\delta}(0) = u_\epsi^0,  \\
		& \dot{u}_{\epsi\delta}(0) = u_\epsi^1, 
	\end{empheq}  
\end{subequations}
 for  some given initial  scaled  lattice  displacement
$u_\epsi^0$  and  velocity $u_\epsi^1$.   Here, the gradient $\partial I_{\epsi\delta}$ of
the smooth energy $I_{\epsi\delta}$ is computed with respect to a
suitable finite dimensional, atomistic, $L^2$-like topology with norm $\| \cdot \|_\epsi$, see \eqref{eq:atomistic inner product} below.  It  corresponds to the system of conservative forces acting  on  the discrete state.    The constant $\nu > 0$  represents   a friction coefficient, whereas $\rho > 0$ denotes the atomic mass. The case $\rho = 0$, corresponding to the purely viscous setting, is covered by our analysis, as well.    In particular, $\nu \dot{u}_{\epsi\delta}$ is a linear friction term, modeling the interaction of the discrete system with the environment.  
System \eqref{eq:intro_gf0} corresponds to the time evolution of an  elastic body under viscous damping. In one and two space dimensions, such  damping
may  arise  from some friction  effect,  see \cite{bonelliAtomisticGrapheneFlakes,  Fundamentals of Friction and Wear on the Nanoscale, krylovAtomisitcMechanismsforFriction} and the references therein. Note that \eqref{eq:intro_gf0} can
be easily complemented with forces, which we however neglect for the sake of
notational simplicity.   
In particular,  the smooth Cauchy problem \eqref{eq:intro_gf0} for ODEs  has a unique
global solution $u_{\epsi\delta}$,  cf.\ Remark \ref{rem:ode}. 

We will prove that (suitable interpolations of)  the atomistic trajectories 
$u_{\epsi\delta}$ converge as $\epsi,\, \delta \to 0$ to the unique  weak 
solution $w \colon  (0,T ) \times  \Omega   \to \mathbb{R}^d$ of 
	\begin{subequations}\label{eq:intro_gf1} 
		\begin{empheq}[left=\empheqlbrace]{align}
				& \rho \partial_{tt} w + \nu \partial_t w - \div (\C : \D^s w) =0 \quad &&\text{in} \
		   (0,T ) \times  \Omega,   \\
		& w = 0  \quad &&\text{on} \
		   (0,T) \times  \partial \Omega,   \\
		& w(0,\cdot)  = w^0 \quad &&\text{in} \ \Omega, \\
		 &   \partial_t w(0,\cdot ) = w^1 \quad &&\text{in} \ \Omega,  
	\end{empheq}
\end{subequations}
where $w^0$  and $w^1$  are   suitable limits of $u^0_\epsi$  and $u^1_\epsi$, respectively. Due to the volume scaling in the energy, $\nu > 0$    is now interpreted as  a friction coefficient, $\rho >0$ is  the  mass density, whereas the case $\rho = 0$ corresponds to a purely viscous flow. 

In order to prove convergence of solutions,  we pass to the limit directly in the atomistic equation  of motion  \eqref{eq:intro_gf0_1}. To this end,  we  reformulate \eqref{eq:intro_gf0_1}  in terms of
an    energy-dissipation-inertia equality  \cite{santambrogioEuclideanMetricWasserstein2016}, namely, 
\begin{align}\label{eq:edie intro}
	 	&\dfrac{\rho}{2} \Vert \dot u_{\epsi\delta} (t) \Vert_\epsi^2  +  I_{\epsi\delta}(u_{\epsi\delta}(t)) +  \dfrac{1}{2\nu} \int_0^t \Vert \rho \ddot u_{\epsi\delta} + \partial I_{\epsi\delta} (u_{\epsi\delta}) \Vert_\epsi^2 \dd s+ \frac{\nu}{2}  \int_0^t	\Vert \dot u_{\epsi\delta}\Vert_\epsi^2 \dd s \nonumber \\
	 	&= 	\dfrac{\rho}{2} \Vert \dot u_{\epsi\delta} (0) \Vert_\epsi^2  +  I_{\epsi\delta}(u_{\epsi\delta}(0)) \qquad \forall t\in [0,T].
\end{align}  
From this we deduce uniform bounds, allowing to extract a weakly converging subsequence by compactness. The main work will then be that of identifying the limit.  The crucial point  is the discussion of the convergence of  the terms $\rho \ddot{u}_{\epsi\delta}$ and $\partial I_{\epsi\delta} (u_{\epsi\delta})$,  where the discrete structure of the atomistic  problem comes into play.  We will see that in the purely viscous case  the analysis is simplified in several  respects.    

For proving energy convergence of solutions, and thus also strong convergence of gradients,  we follow  the classical evolutionary $\Gamma$-convergence approach pioneered by \name{Sandier} and \name{Serfaty}
\cite{ss}.   The energy-dissipation-inertia equation \eqref{eq:edie intro} is passed to the  limit to obtain the  corresponding  
equality   for  the momentum   equation at the continuum level.   This  in turn hinges on the
validity of some  $\liminf$-inequalities, chain rules for the
energies $I_{\epsi\delta}$ and $I$, and the so-called {\it
  well-preparedness} of initial data.  For the validity  of  the chain rule, we require slightly stronger assumptions on the initial data,  compared to the proof of    convergence of solutions.   Taking
advantage of the $\Gamma$\nobreakdash-convergence theory in
\cite{Schmidt_atomisticToCont}, the main point will be that of
checking 
\begin{align}\label{eq:intro_liminf}
  &\int_0^t \|  \rho\partial_{tt} w(s, \cdot)  + \partial I(w(s, \cdot)) \|^2_{L^2(\Omega)}\dd s \nonumber \\
  &\leq \liminf_{\epsi,\delta \to 0} \int_0^t \| \rho \ddot u_{\epsi\delta}(s, \cdot) + \partial
I_{\epsi\delta}(u_{\epsi\delta}(s, \cdot)) \|^2_\epsi \dd s \qquad \forall t \in [0,T],
\end{align} see inequality \eqref{eq:liminf for grad of energy} below.

The energy convergence  entails the   strong convergence of gradients  and prevents  oscillations  in  the sequence of atomistic solutions.  This absence of oscillations can be mechanically interpreted as evidence for the validity of the classical Cauchy-Born hypothesis. In particular, microstructures are not expected to develop.

The paper is organized as follows. In Section~\ref{sec:model}, we
specify the lattice model and state the main result,  namely the convergence of solutions and the evolutionary $\Gamma$\nobreakdash-convergence,  see  Theorem~\ref{thm:main}.  Section~\ref{sec:proof 1} is then  devoted to  the proof of the first part of Theorem~\ref{thm:main},  which is  the convergence of solutions.   We  first  recall  some crucial results
from   \cite{Schmidt_atomisticToCont} in
Section~\ref{subsec:schmidt} and then prove the convergence of solutions  in  Sections~\ref{subsec:grad of energy}--\ref{subsec:conclusion thm 1}. In Section~\ref{sec:proof 2},  we address   the second part of Theorem~\ref{thm:main},  namely energy convergence.  In Section~\ref{subsec:lsc},  we establish the 
 $\liminf$-inequalities,  including \eqref{eq:intro_liminf},    which are then used in Section~\ref{subsec:concluding the proof}  to conclude the argument.

	\section{Setting of the problem and main result}\label{sec:model} 

	In this section, we introduce notation, present preliminaries,
        and state our main convergence result, namely
        Theorem~\ref{thm:main}. At first, we  introduce  the
        atomistic energy.   We follow the setting of  \cite{Schmidt_atomisticToCont},  upon  restricting to the choice of  homogeneous
        Dirichlet boundary conditions  for simplicity.

Throughout the paper, we denote the single contraction among vectors and tensors as $a \cdot b := a_i b_i$,  $(A\cdot b)_{i} = A_{ij} b_j$,  $(A\cdot B)_{ij} := A_{ik} B_{kj}$, $(\C \cdot A)_{ijkl} := \C_{ijkp} A_{pl}$, and $(A\cdot \C)_{ijkl} := A_{ip} \C_{pjkl}$, where we employ \name{Einstein}'s summation convention over repeated indices and tacitly assume that the dimension of the contracting index matches. The double contraction is denoted by  $A:B := A_{ij} B_{ij}$, $(\C:A)_{ij} := \C_{ijkl} A_{kl}$, and $(A:\C)_{ij} := A_{kl} \C_{klij}$ instead, again under the assumption that the dimensions of the contracting indices agree. 

The transposed of the two-tensor $A$ is classically denoted by $A^\top$, whereas we use the notation $\C^t$ to indicate the minor transposition of the four-tensor $\C$, in components $(\C^t)_{ijkl} = \C_{jikl}$.  A four-tensor is called \df{(left) minor symmetric} if $\C^t = \C$, i.e., $\C_{ijkl} = \C_{jikl}$, and \df{major symmetric} if $\C_{ijkl} = \C_{klij}$.

	\subsection{Lattice model}
	
\begin{figure}[t]
		\pgfmathsetmacro{\rows}{7}
		\pgfmathsetmacro{\columns}{11}
		
		\begin{tikzpicture}
			\draw[very thin,color=gray] (0.5,0.5) grid (\columns+0.5,\rows+0.5);
			
			\draw[<->, thin] (1,0.7) --  node[below] {$\epsi$}  (2,0.7) ;
			\draw[<->, thin] (0.7,1) --   node[left] {$\epsi$} (0.7,2);

			\draw (6,4) ellipse (4.9 and 3.2);
			\node[below left] (tildeomega) at (6,7)  {$\ti{\Omega}$};
			
			\begin{scope}[shift={(1,1)}]
				
			\begin{scope}
				\path[draw, thick,use Hobby shortcut,closed=true]
				(1.2,1.2) .. (2.3,1.5) .. (4.5,1.7) .. (6.5,1.3) .. (8.6,3) .. (5.4,4.8) .. (3.3,4.2) .. (1.8,4.4) .. (1.8,2.8);
			\end{scope}
			\node[above right] (omega) at (2,4)  {$\Omega$};
			
			\node[circle,fill=black,inner sep=0pt,minimum size=5pt] (center) at (6.5,3.5) {};
			\draw[->] (center) --+(-0.5,-0.5) node[left, below] {$\epsi z_1$};
			\draw[->] (center) --+(0.5,-0.5) node[right, below] {$\epsi z_2$};
			\draw[->] (center) --+(0.5,0.5) node[right, above] {$\epsi z_3$};
			\draw[->] (center) --+(-0.5,0.5) node[left, above] {$\epsi z_4$};
			
			\draw[fill=black]  (3.5,2.5) circle (1.5pt) node[below] {${\scriptstyle \overline{x}(x,\epsi)}$};
			\draw[fill=black]  (3.8,2.8) circle (1.5pt) node[below, right, xshift = -0.15cm, yshift = -0.15cm] {${\scriptstyle x}$};
			
			\draw[thick, dashed] (3,3) -- (4,3)  node[above, left, yshift= 0.2cm, xshift = 0.1cm]  {${\scriptstyle Q_\epsi(x)}$} ; 
			\draw[thick, dashed] (4,3) -- (4,2);
			\draw[thick] (3,3) -- (3,2) -- (4,2);
			
		\end{scope}
			
		\end{tikzpicture}	
		\caption{Scaled reference configuration of the atomistic system. Here, the lattice is $\L_\epsi = \epsi \Z^2$, i.e., $A$ is the identity matrix. The labeling of the corner points of a scaled cell is indicated by the vectors $\epsi z_i$, $i = 1, \dots, 4$.}
		\label{fig:ref config}
	\end{figure}
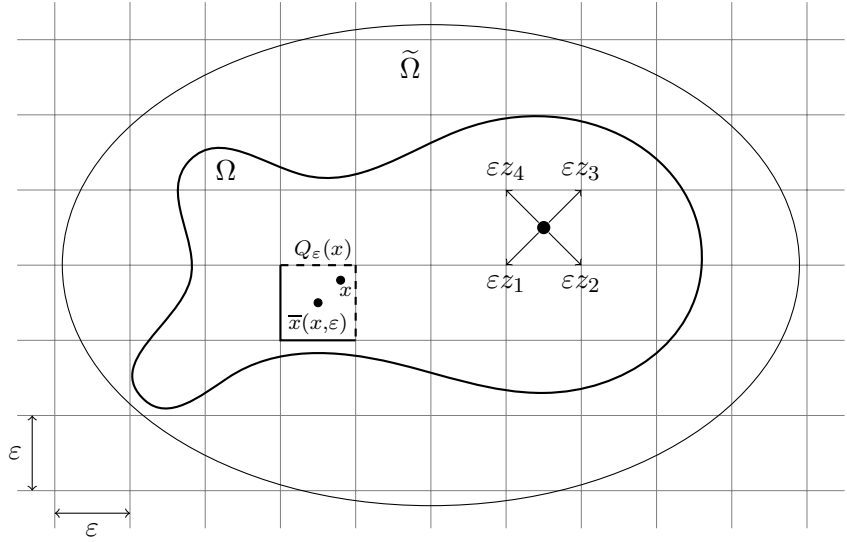

	Let  $v_1, \dots, v_d \subset \R^d$, $ d \geq 1$, be linearly independent vectors.  Define the \emph{Bravais lattice} \[ \L :=  A \cdot \Z^d  := \{\lambda_1 v_1 + \dots + \lambda_d v_d  :  \lambda_1, \dots, \lambda_d \in \Z\}.\] Here, the columns of the matrix $A = (v_1, \dots, v_d)$ consist of the vectors $v_i$, which we assume with no loss of generality to be ordered in such a way that $\det A >0$.  Let the small parameter $\epsi >0$ measure the interatomic distance (see Figure~\ref{fig:ref config}) and define the \df{scaled} Bravais lattice $\L_\epsi  := \epsi A \cdot \Z^d $. An \df{$\epsi$-cell} of the scaled lattice is a set of the form $\epsi A \cdot ( \lambda + [0,1)^d)$ for some $\lambda \in \Z^d$. We denote the set of \df{barycenters of $\epsi$-cells} by $\L'_\epsi := \epsi A \cdot  \left( (1/2, \dots, 1/2) + \Z^d  \right)$. Moreover, we indicate by  $Z = (z_1, \dots, z_{2^d} ) \in \R^{d \times 2^d}$ a fixed \df{labeling} of the $2^d$ corner points of the \df{reference cell} $A\cdot \{-1/2,1/2\}^d$, so that $\L_\epsi = \L'_\epsi + \epsi \{z_1, \dots, z_{2^d}\}$.

	 Let $\Omega \subset \R^d$  be a bounded domain  with $C^{1,1}$ boundary.   For every $x \in \Omega$, we denote by $\overline{x}(x, \epsi) \in \L'_\epsi$ the barycenter of the $\epsi$-cell containing  $x$ and correspondingly by ${Q_\epsi(x) := \overline{x}(x,\epsi) +  A \cdot [- \epsi/2, \epsi/2 )^d}$ the $\epsi$-cell containing $x$, see Figure~\ref{fig:ref config}.  In particular, for $\bar{x} \in \mathcal{L}'_\eps$, we have $Q_\epsi(\ove{x}) =  \overline{x} +  A \cdot [- \epsi/2, \epsi/2 )^d$ and for $x \in \mathcal{L}_\eps$, we have $Q_\epsi({x}) = x + A \cdot [0,\epsi)^d$.

 	At the atomistic level, the central  objects  are \emph{lattice deformations} $y \colon \L_\epsi \cap \Omega \to \R^d$, describing  the  response of each lattice point in $\Omega$ to mechanical actions.  We  impose homogeneous Dirichlet conditions on the whole boundary  in the following sense:  let $\ti{\Omega} \subset \R^d$   be  a    bounded open set with $\Omega \subset \subset \ti{\Omega}$  (compactly contained). Given a lattice deformation $y \colon \L_\epsi \cap \Omega \to \R^d$, we  extend it   to $\L_\epsi \cap \ti{\Omega}$ by setting $y(x) = x$ for $x \in \L_\epsi \cap (\ti{\Omega} \setminus \Omega)$. This choice of boundary conditions is meant to simplify notations.  It  is  possible to also consider nonhomogeneous Dirichlet conditions or even Neumann-like conditions, see \cite{Schmidt_atomisticToCont}. This would  however  require some additional modeling choice and a more technical treatment. 

	We denote the set of  \df{barycenters of $\epsi$-cells inside $\ti{\Omega}$} by $\L'_\epsi (\ti{\Omega}) := \{ \overline{x} \in  \mathcal{L}_\epsi'  \colon  \ove{Q_\epsi (\ove{x})} \subset \ti{\Omega} \}$  and abbreviate the union of cells which are fully contained in $\ti{\Omega}$  by $\bigcup Q_{\epsi} := \bigcup_{\ove{x} \in \L'_{\epsi}(\ti{\Omega})} Q_{\epsi}(\ove{x})$.  For $y\colon \L_\epsi\cap \ti{\Omega} \to \R^d$ and  $x \in Q_\epsi(\overline{x})$ with $\overline{x} \in \L'_\epsi (\ti{\Omega})$,  we  collect  the deformations  $y_i = y_i (x) := y(  \overline{x}  + \epsi z_i)$  of corner points,  and define the   \df{discrete gradient} $\bar{\D} y \colon   \ti{\Omega}  \to \R^{d \times 2^d}$ of $y$  as
		\begin{align}\label{eq:def of disc grad}
			\bar{\D} y(x) := \epsi^{-1} (y_1 - \overline{y}, \dots, y_{2^d} - \overline{y}) \quad \mathrm{with}\quad   \overline{y} := \frac{1}{2^d} \sum_{i = 1}^{2^d} y_i,
		\end{align}
		for $x \in \bigcup Q_{\epsi}$ and $\bar{\D} y(x) = 0$ for $x \in \ti{\Omega}\setminus    \bigcup Q_{\epsi}$.  	As the discrete gradient is constant on individual cells, in the following we often indicate its argument simply by $\ove{x}$,  that is, the barycenter of such  a  cell. 
	
	The   energy is defined  in terms of the \df{lattice displacement}, namely, $u := y - \id$, where   $\id$ indicates the identity mapping on $\R^d$.    In order to address    the regime of small deformations,  we consider the \df{scaled lattice displacement} $u = \delta^{-1} (y - \id)$, where $\delta >0$ is a small parameter, and define the set of \df{admissible (scaled lattice) displacements} as 
	\begin{align*}
		\A_\epsi := \big\{ u: \L_\epsi \cap \ti{\Omega}  \to \R^d : u(x) = 0 \text{  for  } x \notin \Omega  \big\}.
	\end{align*}
 Then,  the \df{scaled atomistic energy} $I_{ \epsi \delta } : \A_\epsi \to \R$ is defined as
	\begin{align*}
		I_{\epsi\delta} (u) := \frac{\epsi^d  \det A }{\delta^2} \sum_{\ove{x} \in \L'_\epsi (\ti{\Omega}) } W(Z + \delta \bar{\D} u(\ove{x})), 
	\end{align*}
	where the \df{elastic cell energy} $W$ of the configuration
        only depends on the discrete gradient. The factor $\epsi^d
        \det A$ in the definition of $I_{\epsi\delta}$ corresponds to
        the volume of a single scaled cell, whereas the scaling
        $\delta^{-2}$ is instrumental  for the setting of small
        displacements.   	Let $\SO(d) := \{ Q \in \R^{d
          \times d}  \colon   Q^\top \cdot Q  =  {\rm Id},  \, \det Q = 1  \}$ be the standard special orthogonal group. 
      
	\newpage
	\begin{assump} \label{assumptions}\leavevmode 	We impose the following assumptions on the  cell  energy $W$. 
		\begin{enumerate}[label=\emph{(\roman*)}] 
			\item The  cell  energy $W\colon \R^{d \times 2^d} \to [0, \infty]$ is  invariant under  translations and rotations, that is, given $F \in \R^{d \times 2^d}$ we have that  
			\begin{align*}
				W \big(R\cdot F +  \underbrace{(c,\dots, c)}_{2^d\mathrm{-times}}\big) = W(F)
			\end{align*} holds for all $R \in \SO(d)$ and $c \in \R^d$.
			\item $W( F  ) =0$ if and only if there exists a rotation $R \in\SO(d)$ and a translation $c\in \R^d$ such that $ F_i  = R\cdot z_i +c $. 			
			\item $W$ is $\Cc^2$ and the Hessian $\mathbb{H} := D^2 W(Z) \in \R^{d \times 2^d \times d \times 2^d}$ is positive definite on the orthogonal complement of the subspace spanned by infinitesimal translations $(x_1, \dots, x_{2^d}) \mapsto (c, \dots, c)$  and infinitesimal rotations $(x_1, \dots, x_{2^d}) \mapsto (S\cdot x_1, \dots, S\cdot x_{2^d})$, where $c \in \R^d$ and $S \in \R^{d\times d}_\mathrm{skew} := \{   F  \in \R^{d \times d}  \colon F + F^\top = 0\}$. 
			\item $W$ grows at infinity at least quadratically on the orthogonal complement of the subspace spanned by infinitesimal translations, i.e., 
			\begin{align*}
				\liminf_{\substack{\vert F \vert \to \infty, \\ F \in V} } \frac{W(F)}{\vert F \vert^2} > 0,
			\end{align*}	
			with $V = \{ F \in \R^{d \times 2^d} \colon  F_1 + \dots + F_{2^d} = 0  \}$.
			\item The gradient of the cell energy grows   at most quadratically, i.e., $ \vert D W( F  ) \vert \leq c\vert F \vert^2 + c$  for all $F \in \R^{d \times 2^d}$  for  some  constant $c > 0$.		
		\end{enumerate}
	\end{assump}	
	
 	Assumptions~\ref{assumptions} (i)--(iv) have already been considered by \name{Schmidt} \cite{Schmidt_atomisticToCont} in the stationary setting and allow to characterize the $\Gamma$-limit of $I_{\epsi\delta}$ as $\epsi,\,\delta \to 0$. In this paper, we strengthen the setting  by additionally assuming Assumption~\ref{assumptions} (v).  This  is  needed for   dealing with  the  gradients of the energies. 
 	
	\subsection{Main result}\label{sec:main result}
	 We start by defining the \df{atomistic $L^2$-inner product} $( \cdot , \cdot )_\epsi$ on the space of scaled lattice displacements. Given $u_\epsi, v_\epsi \colon \L_{\epsi} \cap \ti{\Omega} \to \R^d$, we define \begin{align}\label{eq:atomistic inner product} 
		(u_\epsi, v_\epsi)_\epsi := \epsi^d \det A \sum_{ x \in  \L_{\epsi} \cap \ti{\Omega}} u_\epsi (x) \cdot v_\epsi (x)
	\end{align}
	and denote the  corresponding \df{atomistic $L^2$-norm} by   $\Vert u_\epsi \Vert_\epsi := \sqrt{(u_\epsi , u_\epsi)_\epsi}$.

  We look for a solution  $ u_{\epsi\delta}  \colon [0,T] \times (\L_\epsi \cap \ti{\Omega}) \to \R^d$ with $ u_{\epsi\delta}   (\cdot, x) \in \Cc^2([0,T]; \R^d)$ for all $x \in \L_\epsi \cap \ti{\Omega}$, and $ u_{\epsi\delta}   (t, \cdot) \in \A_\epsi$ for all $t \in [0,T]$ to the \df{atomistic equation of motion}          
	\begin{subequations}\label{eq:atomistic the second}
	\begin{empheq}[left=\empheqlbrace]{align}
			& \rho \ddot{u}_{\epsi\delta}(t, x) + \nu \dot{u}_{\epsi\delta}(t,x )  + \partial I_{\epsi\delta} ( u_{\epsi\delta}   (t,x))= 0 \quad &&\forall x \in \L_\epsi \cap \ti{\Omega}, \ \forall t \in (0,T),  \\
			& u_{\epsi\delta}   (0,x) = u_\epsi^0(x) \quad &&\forall x \in \L_\epsi \cap \ti{\Omega},\\
			& \dot{u}_{\epsi\delta}   (0,x) = u_\epsi^1(x) \quad &&\forall x \in \L_\epsi \cap \ti{\Omega},
	\end{empheq} 
\end{subequations}
where  $u_\epsi^0,  u_\epsi^1  \in \A_\epsi$  are  given \df{initial scaled lattice displacement} and  \df{velocity}, respectively,   $\rho > 0$ corresponds to the atomic mass, and $\nu > 0$ is  a friction coefficient.  The case $\rho= 0$ is also considered and corresponds to the purely viscous setting. 
 Here, $\partial I_{\epsi\delta}( u_{\epsi\delta}  ) $ indicates the first variation of $I_{\epsi\delta}$ in $\A_\epsi$ at $ u_{\epsi\delta}  $,  that is,
 \begin{align*}
 	&\partial I_{\epsi\delta} : (\A_\epsi, \Vert \cdot \Vert_\epsi ) \to (\A_\epsi, \Vert \cdot \Vert_\epsi )^* \nonumber \\
 	&u_{ \epsi\delta}\mapsto \partial I_{\epsi\delta} (u_{\epsi\delta}),
 \end{align*}
where,  for each $v_\epsi \in \A_\epsi$,  
\begin{align}\label{eq:defin of gradient}
	\langle\partial I_{\epsi\delta} (u_{\epsi\delta}), v_\epsi \rangle :=& \frac{\dd}{\dd h} I_{\epsi\delta}( u_{\epsi\delta}   + h v_\epsi) \big\vert_{h = 0} \nonumber \\
	=&  \epsi^d  \det A  \sum_{\ove{x} \in  \L'_{\epsi} (\ti{\Omega})} \frac{1}{\delta} DW(Z + \delta\bar{\D} u_{\epsi\delta}(\ove{x})) : \bar{\D} v_\epsi(\ove{x}).
\end{align}
 The  homogeneous Dirichlet boundary conditions  entail the validity of a discrete divergence theorem, see Lemma~\ref{lem:summation by parts}.
Hence, we can also interpret $\partial I_{\epsi\delta} (u_{\epsi\delta}): \L_\epsi \cap \ti{\Omega} \to \R^d$ as the unique scaled lattice  strain  satisfying
\begin{align}\label{eq:defin disc grad as lattice disp}
		\langle\partial I_{\epsi\delta} (u_{\epsi\delta}), v_\epsi \rangle  = 	(\partial I_{\epsi\delta} (u_{\epsi\delta}), v_\epsi )_\epsi \quad \quad \text{for all $v_\epsi \in \A_\epsi$.}
\end{align}   
As $\L_\epsi \cap \ti{\Omega}$ is a finite point set, problem \eqref{eq:atomistic the second} is nothing but a Cauchy problem for a smooth ODE system,  see  also Remark~\ref{rem:ode} below. It hence admits a unique  global  solution. Indeed, local well-posedness of \eqref{eq:atomistic the second} follows from the classical Picard-Lindelöf theorem.  From the energy-dissipation-inertia equality \eqref{eq:edie intro} corresponding to the atomistic problem \eqref{eq:atomistic the second} one can derive bounds ensuring   that $I_{\epsi\delta}( u_{\epsi\delta}  )$ remains bounded along the trajectories. From the coercivity in Assumption~\ref{assumptions}~(iv) one obtains that trajectories are bounded,    see Theorem~\ref{thm:gamma convergence and compactness} below.  Thus, the solution is  global in time.

	\begin{remark}\label{rem:ode}
			To clarify that the atomistic  equation of motion 
                         \eqref{eq:atomistic the second} can
                        indeed be interpreted as  an  ODE system, let
                        $N_\epsi$ be the number of atoms in $\L_\epsi
                        \cap \Omega$ and $x_1, \dots, x_{N_\epsi}$
                       some  corresponding labeling.  Given a scaled
                        lattice displacement $ u_{\epsi\delta}
                          \in \A_\epsi$, let $ U_{\epsi\delta}
                        := (u_1, \dots, u_{N_\epsi}) \in \R^{d \times N_\epsi}$ be the vector of the displaced atoms, where $u_i =  u_{\epsi\delta}   (x_i)$ for $i = 1, \dots, N_{\epsi}$. This defines a linear isomorphism $\iota \colon \A_\epsi \to \iota(\A_\epsi) \subset \R^{d \times N_\epsi}$. 
			
			The scaled atomistic energy can be
                        equivalently rewritten as
                        $\ti{I}_{\epsi\delta} ( U_{\epsi\delta} ) := I_{\epsi\delta}(\iota^{-1} ( U_{\epsi\delta} ) ) = I_{\epsi\delta} ( u_{\epsi\delta} )$. 			Then, problem~\eqref{eq:atomistic the second} is equivalent to finding $ U_{\epsi\delta}  (t) = \iota( u_{\epsi\delta}  (t))$ satisfying  
		\begin{equation*}
			\left\lbrace \begin{aligned}
				&\rho \ddot{U}_{\epsi\delta} (t)+ \nu \dot{U}_{\epsi\delta}  (t) + \D \ti{I}_{\epsi\delta} ( U_{\epsi\delta}  (t)) = 0 &&  \forall t \in (0,T), \\
				& U_{\epsi\delta}  (0) = U_\epsi^0,\\
				& \dot{U}_{\epsi\delta}  (0) = U_\epsi^1,
			\end{aligned} \right. 
		\end{equation*}
		where $U_\epsi^0 = \iota(u_\epsi^0)$,  $U_\epsi^1 = \iota(u_\epsi^1)$,  and  $\D\ti{I}_{\epsi\delta}  =   \iota(\partial I_{\epsi\delta}( u_{\epsi\delta} ))$.  		
		\end{remark}
	
	The goal of this paper is to show that solutions $
        u_{\epsi\delta} $ to the  atomistic equation of motion~\eqref{eq:atomistic the second}  converge  as
        $\epsi, \, \delta \to 0 $  to  the unique weak  solution of the  \df{linear continuum momentum balance}  
	\begin{subequations}\label{eq:continuous problem}
		\begin{empheq}[left=\empheqlbrace]{align}
			&  \rho \partial_{tt} w + \nu \partial_t  w  - \div (\C : \D^s w)=0    &&\text{in} \
			(0,T ) \times  \Omega,    \\
			& w  = 0  &&\text{a.e. in  } (0,T) \times \partial \Omega,\\
			& w  (0,\cdot) =  w^0   &&\text{a.e. in  } \Omega,\\
			& \partial_t w  (0,\cdot) =  w^1  &&\text{a.e. in  } \Omega,
		\end{empheq}
	\end{subequations}
 where    $ w^0  \in  H^1_0  (\Omega;\R^d)$ and $w^1  \in L^2(\Omega;\R^d)$   are the given  \df{initial continuum displacement}   and the \df{initial continuum velocity}, respectively.   The fourth-order symmetric  \emph{elasticity tensor} $\C$ is  defined  as
	\begin{align}\label{eq:tensor}
		\C := Z \cdot  \Hz^t \cdot Z^\top =  Z \cdot D^2 W(Z)^t \cdot Z^\top \in \R^{d \times d \times d \times d},
	\end{align} and  we recall that  $\D^s  w  := (\D  w  + \D 
	w^\top)/2$ denotes the symmetrized gradient. In
	particular, $\C$ is solely determined by microscopic
	quantities and $\C$ is minor and major symmetric, see
	Assumptions~\ref{assumptions}~(iii) and
	Lemma~\ref{lem:symmetry of tensor} below. 

   We say  that $w \in L^2(0,T;H^1_0(\Omega;\R^d)) \cap  H^1(0,T;L^2(\Omega;\R^d)) \cap H^2(0,T;H^{-1}(\Omega;\R^d))$ is a \df{weak solution to \eqref{eq:continuous problem}} if it satisfies   $w(0) = w^0$, $\partial_t w (0) = w^1$, and  
	\begin{align*}
		\la \rho \partial_{tt} w , v \ra_{H^1_0} + (\nu \partial_t w , v )_{L^2} + \int_\Omega \D^s w : \C : \D^s v \dd x = 0 \quad \text{a.e.\ in $(0,T)$},
	\end{align*} 
	and for all $v \in H^1_0(\Omega;\R^d)$,  where  $( \cdot , \cdot)_{L^2}$ denotes the $L^2$-inner product	 and  $\la \cdot, \cdot \ra_{H^1_0}$ the duality in $H^1_0$.  	In system \eqref{eq:continuous problem}, the term $- \div(\C : \D^s w)$   is hence   interpreted in the weak  $H^1$-sense and corresponds to the  $H^1$-gradient $\partial I$ of $I$, i.e.,  the 
        Fr\'echet differential     of  the  \df{linearized continuum energy} $I \colon H_0^1(\Omega; \R^d) \to \R$ given by 
	\begin{align}\label{eq:cont deformation energy}
		 I( w )&:= \frac{1}{2} \int_\Omega (\D  w  \cdot Z):\mathbb{H} :(\D  w  \cdot Z)  \dd x
		 = \frac{1}{2} \int_\Omega \D^s  w   :\C : \D^s  w  \dd x.
	\end{align} 
 In particular,  $ \partial I : H_0^1(\Omega;\R^d) \to H^{-1}(\Omega;\R^d)$ is given by
\begin{align*}
	\langle \partial I(w), v \rangle = \frac{\dd }{\dd h} I(w + hv) \big\vert_{h = 0} = \int_\Omega (\D w  \cdot Z) : \Hz : (\D v \cdot Z)  \dd x  = \int_\Omega \D^s w : \C : \D^s v \dd x. 
\end{align*}
We derive the formula above   in Section~\ref{subsec:grad of energy} below, see \eqref{eq: forrrlater}.

	In order to formulate our statement, we need to specify a notion of convergence
        for a sequence of scaled lattice displacements $u_\epsi\colon
        \L_{\epsi} \cap \ti{\Omega} \to \R^d$. Let $ w  \in L^2(\Omega; \R^d)$ be a   continuum  displacement and define its local mean on the lattice cells by
	\begin{align}\label{eq:definition P}
		P_\epsi  w  (x) := \fint_{Q_\epsi(x) } \ti{ w  }(\xi) \dd \xi \quad \forall x \in \L_\epsi \cap \ti{\Omega},
	\end{align}
	where $\ti{ w } \colon   \ti{\Omega}  \to \R^d$ is the trivial extension of $ w $ to $ \ti{\Omega} $. We say that $u_\epsi$  \df{\textup{AC}-converges} to $ w $   and write  $u_\epsi \acto  w $  (\textup{AC} stands for \textit{atomistic-to-continuum}  convergence)  as $\epsi \to 0$ if 
	\begin{align}\label{eq:convergence}
		\Vert u_\epsi - P_\epsi  w   \Vert_\epsi \to 0. 
	\end{align}

	The main result of this paper  is  the following.

	\begin{theorem}[Discrete-to-continuum linearization]\label{thm:main}
	 	Suppose that Assumption~\emph{\ref{assumptions}}   holds and  let $\epsi_k$ and $\delta_k$ be such that  $\epsi_k \to
	 	0$ and  $\delta_k \to 0$.

	 	\noindent\underline{\emph{Convergence of solutions:}}  Let	$u_\epsi^0,   u_\epsi^1  \in \A_\epsi$ be an initial scaled lattice
	displacement  and velocity,   and   let  $ w^0 \in
	 H^1_0 (\Omega;\R^d)$, $w^1  \in L^2(\Omega;\R^d)$  be an initial continuum displacement  and velocity,    respectively.    If   $ u_{\epsi_k}^0  \acto  w^0 $,
	 $ u_{\epsi_k}^1  \acto  w^1 $ as  $k \to \infty$   and
	$\sup_k   I_{\eps_k\delta_k}  (u_{\epsi_k}^0) < + \infty$, then the solutions
	$( u_{\epsi_k\delta_k} )_k$ to the  atomistic equation of motion  \eqref{eq:atomistic the second}  \textup{AC}-converge  to $ w $ pointwise in time, where $ w $ is the unique   weak  solution to the     continuum momentum equation  \eqref{eq:continuous problem}.

	\noindent\underline{\emph{Energy convergence:}}  If, additionally, $ w^0 \in	H^2(\Omega;\R^d)$, $w^1  \in H^1(\Omega;\R^d)$,  and 
	$I_{ \epsi_k \delta_k}( u_{\epsi_k}^0 ) \to I( w^0 )$ as  $k \to \infty$, then the solutions
	$( u_{\epsi_k\delta_k} )_k$ to the  atomistic equation of
        motion    additionally satisfy  the convergences   $
         I_{\eps_k\delta_k}  ( u_{\epsi_k\delta_k}  (t)) \to
        I(w(t))$ for all $t \in [0,T]$ and $\bar{\D} 
        u_{\epsi_k\delta_k}  (t) \to \D w(t)\cdot Z$ strongly in $L^2(\Omega;\R^d)$ for  all  $t \in (0,T)$.   Moreover, $w$ is a strong solution of  \eqref{eq:continuous problem}. 
	\end{theorem}

	In order to prove the  first part of Theorem~\ref{thm:main}, i.e., the convergence of solutions,  we pass to the limit in the atomistic equation of motion  \eqref{eq:atomistic the second}. To this end, we  need to  deduce uniform bounds on the individual terms of the atomistic equation of motion.  This is achieved by  the energy-dissipation-inertia equality \eqref{eq:edie intro},  which in the 
	notation of Theorem \ref{thm:main} reads as 
     \begin{align}
  &\frac{\rho}{2}\|\dot u_{\epsi_k\delta_k}(t)\|^2_{\epsi_k} + I_{\epsi_k\delta_k}(u_{\epsi_k\delta_k}(t)) + \frac{\nu}{2}\int_0^t   \|\dot u_{\epsi_k\delta_k}(s)\|^2_{\epsi_k} \, \dd s  \nonumber \\
  &\quad + 
    \frac{1}{2\nu} \int_0^t \| \rho \ddot u_{\epsi_k\delta_k} (s) +
    \partial I_{\epsi_k\delta_k}(u_{\epsi_k\delta_k}(s))\|^2_{\epsi_k}  \, \dd s =  \frac{\rho}{2}\|u^1_{\epsi_k}\|^2_{\epsi_k} + I_{\epsi_k\delta_k}(u^0_{\epsi_k})  \label{eq:zoom}
     \end{align}
     for all $t \in [0,T]$, see Section \ref{sec:apriori}  for its derivation.

 On the contrary, we  prove  the second part of  Theorem~\ref{thm:main} by following the classical
        evolutive  $\Gamma$\nobreakdash-convergence approach pioneered by
        \name{Sandier} and \name{Serfaty} \cite{ss},  see also
         \cite{Mielke,
       Ortner_GradFlowAsSelection}.
      This hinges on  the fact that the atomistic evolution problem
     \eqref{eq:atomistic the second}  is actually  equivalent to \eqref{eq:zoom}.  One passes to the $\liminf$
     in the latter as $k \to \infty$. The right-hand side converges by
     assumption to $\frac{\rho}{2}\|w^1\|^2_{L^2(\Omega)} + I
     (w^0)$. As regards the terms in the left-hand side, we will
     check  that, for each $t \in [0,T]$, 
     \begin{subequations}\label{eq:zoom1}
       \begin{align}
         & \|\partial_t w(t)\|^2_{L^2(\Omega)} \leq \liminf_{k \to
                                                           \infty} \|\dot u_{\epsi_k\delta_k}(t)\|^2_{\epsi_k}, \label{eq:liminf for time derivative}\\[1.5mm]
          & \int_0^t   \|\partial_t w (s)\|^2_{L^2(\Omega)} \, \dd s\leq \liminf_{k  \to
                                                                                \infty}  \int_0^t   \|\dot u_{\epsi_k\delta_k}(s)\|^2_{\epsi_k} \, \dd s, \label{eq:liminf for integral of time derivative}\\[1.5mm]
                                                                                        &I(w(t)) \leq \liminf_{k \to
                   \infty} I_{\epsi_k\delta_k}(u_{\epsi_k\delta_k}(t)) , \label{eq:zoom12}\\
         &  \int_0^t \| \rho \partial_{tt} w (s) +
           \partial I (w(s))\|_{L^2(\Omega)}^2  \, \dd s \notag \\ & \quad 
          \leq \liminf_{k \to  \infty}
                                                         \int_0^t \| \rho \ddot u_{\epsi_k\delta_k} (s) +
                                                        \partial
           I_{\epsi_k\delta_k}(u_{\epsi_k\delta_k}(s))\|_{\epsi_k}^2
           \, \dd s. \label{eq:liminf for grad of energy}
       \end{align}
     \end{subequations} 
Then, taking the $\liminf$ as $k \to \infty$ in equality \eqref{eq:zoom} will entail that the following
energy-dissipation-inertia (in)equality holds
\begin{align*}
&\frac{\rho}{2}\|  \partial_t   w(t)\|^2_{L^2(\Omega)} + I(w(t))+ \frac{\nu}{2}\int_0^t   \|\partial_t w (s)\|^2_{L^2(\Omega)} \, \dd s
 \nonumber \\
  &\quad +  \frac{1}{2\nu} \int_0^t \| \rho \partial_{tt} w (s) +
  \partial I (w(s))\|_{L^2(\Omega)}^2  \, \dd s   \leq   \frac{\rho}{2}\|w^1\|^2_{L^2(\Omega)} + I
     (w^0) \qquad \forall t \in [0,T].\label{eq:zoom2}
\end{align*}
As long as the initial conditions $w(0)=w^0$ and $\partial_t w(0)=w^1$
are satisfied, the latter is indeed equivalent to the corresponding equality, as well
as to linear continuum momentum balance \eqref{eq:continuous
  problem}.  We refer to Section \ref{sec:apriori} below where the precise argument is given for the discrete problem.  This additionally implies that each $\liminf$ in
\eqref{eq:zoom1} is actually a limit. In particular, the energy
$I_{\epsi_k\delta_k}(u_{\epsi_k\delta_k}(t))$ converges to $I(w(t))$ at
all times, which eventually implies that $\bar{\D} u_{\epsi_k\delta_k}(t) \to \D
w(t)\cdot Z$ strongly in $L^2(\Omega;\R^{d\times 2^d})$, as well.

	\begin{remark}[Purely viscous case]
		Theorem~\ref{thm:main} above directly   extends  to the purely viscous setting of $\rho = 0$. In this case, the continuum problem \eqref{eq:continuous problem} corresponds to a gradient flow  and, assuming the milder regularity assumptions $w^0 \in  H^1_0 (\Omega;\R^d)$ and $w^1 \in L^2(\Omega;\R^d)$,    we can show that $\div (\C : \nabla^s w ) \in L^2(0,T;L^2(\Omega;\R^d))$, see Remark~\ref{rem:purely viscous setting strong solutions} below. Hence, $-\div (\C : \D^s w)$ is the $L^2$-gradient of the elastic continuum energy $I$, given in \eqref{eq:cont deformation energy}. In particular, we obtain strong solutions for the continuum problem and strong convergence of gradients without the  higher regularity assumption on the initial data.    The proof can be simplified  as inequality \eqref{eq:liminf for grad of energy} directly follows from Lemma~\ref{lem:convergence} below.   We refer to Remarks \ref{rem:inertia grad in L2}  and \ref{rem:purely viscous setting strong solutions} for details.  
	
	\end{remark}

	\section{Convergence of solutions}\label{sec:proof 1}
		In this section, we prove the  first part of  Theorem~\ref{thm:main},  namely, the convergence of solutions.  To this end, we derive uniform bounds from the energy-dissipation-inertia equality  \eqref{eq:zoom}  and use compactness to obtain converging subsequences. The main  focus  is then on identifying limits.
		
		For the reader's convenience, the argument is divided into steps. 
		In Section~\ref{subsec:schmidt}, we start  by recalling the compactness and  $\Gamma$\nobreakdash-convergence result of \cite{Schmidt_atomisticToCont}. 	In Section~\ref{subsec:grad of energy},  $\partial I$ is computed and we further peruse properties of scaled lattice displacements.   In Section~\ref{sec:apriori}, we derive \eqref{eq:zoom} and state the resulting a priori estimates. 	Sections~\ref{subsec:limit2} and \ref{subsec:limit1} are then devoted to computing the limits of  $\nu \dot{u}_{\epsi_k\delta_k}$ and $\rho \ddot{u}_{\epsi_k\delta_k} + \partial I_{\epsi_k\delta_k}(u_{\epsi_k\delta_k})$,  respectively. We conclude the proof in Section~\ref{subsec:conclusion thm 1}.

		In the remainder of the paper, given the
		sequences  $\epsi_k, \delta_k \to 0$ with $k \to
		\infty$   from Theorem~\ref{thm:main}, we   use the short-hand
		notation $I_k :=  I_{\epsi_k \delta_k}$, 
		$u_k :=  u_{\epsi_k\delta_k}$,  $(\cdot , \cdot)_{k} := ( \cdot, \cdot)_{\epsi_k}$, and $\Vert \cdot \Vert_k := \Vert \cdot \Vert_{\epsi_k}$. Throughout this  paper,  $c$ denotes a positive constant, independent of $k$, but possibly varying from line to line.   
		 Moreover, we will denote the standard $L^2$-inner product by $( \cdot , \cdot)$ and denote duality pairings between the space $X$ and its dual $X^*$ by $\la \cdot , \cdot \ra_X : X^* \times X \to \R$. For $X = H^1_0(\Omega;\R^d)$, we will simply write $\la \cdot, \cdot \ra_{H^1_0} =: \la \cdot, \cdot \ra.$

			\subsection{$\Gamma$-convergence and  compactness}\label{subsec:schmidt}
		
		In this section, we briefly summarize the $\Gamma$\nobreakdash-convergence and compactness result  by  \name{Schmidt} \cite{Schmidt_atomisticToCont}, for these  are crucial ingredients in the proof of Theorem~\ref{thm:main}. 
		\begin{theorem}{\cite[Theorem 2.6]{Schmidt_atomisticToCont}} \label{thm:gamma convergence and compactness}
			Under Assumptions~\emph{\ref{assumptions}} we have the following: 
			
			\smallskip
			\noindent\underline{\emph{Compactness:}} If $I_k(u_k)$ is equibounded for some sequence  $(u_k)_k$ of  scaled lattice displacements, then there exists a (not relabeled) subsequence such that $u_k \acto  w $ for some $ w  \in H_0^1(\Omega; \R^d)$. Moreover,   $\bar{\D} u_k \wto \D \ti{ w } \cdot Z$ weakly in $L^2( \ti{\Omega} ; \R^{d \times 2^d})$,  where $\ti{ w } \colon  \ti{\Omega}  \to \R^d$ is the trivial extension of $ w $ to $\ti{\Omega}$.

			\smallskip
			
			\noindent\underline{\emph{$\Gamma$-convergence:}} The functionals $I_k$ $\Gamma$-converge to the functional $I$  in the sense of \textup{AC}-convergence, i.e., the following  is  satisfied:
			\begin{enumerate}[label=\emph{(\roman*)}] 
				\item \emph{$\liminf$-inequality:} Every sequence $(u_k)_k$,  $u_k \in  \mathcal{A}_{\epsi_k}$,   \textup{AC}-converging  to some $ w  \in H_0^1( \Omega;\R^d)$ satisfies the estimate \[  \liminf_{k \to \infty} I_k (u_k) \geq I( w ).  \]  
				\item \emph{Recovery sequence:}  For every $ w  \in H_0^1( \Omega; \R^d)$ there exists a sequence $(u_k)_k$, $u_k \in \mathcal{A}_{\epsi_k}$, such that $u_k \acto  w  $ and \[ \lim_{k \to \infty} I_k(u_k) = I( w ). \]
			\end{enumerate}
		\end{theorem}

		\begin{remark}\label{rem:bdry and Cauchy Born}
			The $\Gamma$-convergence and compactness result  can be extended in various directions: 

			a) For recovery sequences, that are sequences $u_k \acto  w  \in H_0^1(\Omega; \R^d)$ such that $I_{k}(u_k) \to I( w )$, one can also obtain strong convergence of discrete gradients, see \cite[Theorem~2.7]{Schmidt_atomisticToCont}. 

b) AC-convergence  \eqref{eq:convergence} is equivalent to the classical $L^2$-convergence for suitable interpolations,  see for instance  Lemma~\ref{lem:conv of norms} below. The $\Gamma$-convergence result above could be equivalently reformulated to hold with respect to the strong $L^2$-  or  the weak $H^1$-topology for such interpolations. 

c) The compactness statement in Theorem~\ref{thm:gamma convergence and compactness} hinges on the fact that Dirichlet boundary conditions are imposed on the whole $\partial\Omega$. By imposing Dirichlet conditions on a part of $\partial\Omega$ only, one could only conclude local weak $L^2_\loc$-convergence of the discrete gradients. As already mentioned, other boundary conditions could be considered as well, by resorting to the theory in \cite{Schmidt_atomisticToCont}.
		\end{remark}

\subsection{Gradients of the energy and auxiliary results}\label{subsec:grad of energy}

		 In order to prove convergence of solutions, we aim at  showing the weak convergence of the Fr\'echet differential $\partial I$  in  a suitable  dual space. Due to the varying underlying spaces,   implementing  this strategy  requires some care. We  show that, for a suitable sequence $v_k \acto v$, the G\^{a}teaux derivatives in direction $v_k$ converge (Lemma~\ref{lem:convergence} below).  
		To this end, we introduce an approximation of test functions by scaled lattice displacements and study their  properties.

		We start by computing the Fr\'echet differential  of $I$  in Section~\ref{subsec:diff of I}. In Section~\ref{subsec:scaled lattice displacements}, we  address the approximation by scaled lattice displacements.

		\subsubsection{Fr\'echet differential of $I$}\label{subsec:diff of I}	
		Recall that the Hessian of $W$ evaluated at $Z$ is $\mathbb{H} := D^2 W(Z)$ and therefore $\Hz$ is major symmetric, namely, $\Hz_{ijkl} = \Hz_{klij}$.
		
		Given $ w  \in H_0^1(\Omega; \R^d)$ and $v \in \Cc^\infty_\comp(\Omega; \R^d)$, the variation of $I$  in  $H^1(\Omega;\R^d)$  at $ w $ in direction $v$ can be classically computed as
		\begin{align*}
			\la \partial I( w ) , v \ra &= \dfrac{\dd}{\dd h} I( w  + hv)  \big\vert_{h = 0} \\
			&=  \frac{1}{2} \int_\Omega \dfrac{\dd}{\dd h}\Bigl( \bigl( (\D  w  + h \D v)  \cdot Z  \bigl) :  \Hz  : \bigl( (\D  w  + h \D v)  \cdot Z \bigr) \Bigr) \Big\vert_{h = 0} \dd x.
		\end{align*}
		Due to the major symmetry of $\Hz$, one readily checks that
		\begin{align*}
			\dfrac{\dd}{\dd h}\Bigl( \bigl( (\D  w  + h \D v)  \cdot Z  \bigl) :  \Hz  : \bigl( (\D  w  + h \D v)  \cdot Z \bigr) \Bigr) \Big\vert_{h = 0}=  2  (\D  w  \cdot Z) : \mathbb{H} :(\D v \cdot Z).
		\end{align*}
		Therefore, we deduce that 
		\begin{equation}\label{eq:derivative of limiting energy in terms of H}
			\begin{aligned}
				\la \partial I( w ) , v \ra =  \int_\Omega   (\D  w  \cdot Z) :\mathbb{H}: (\D v \cdot Z) \dd x.
			\end{aligned}
		\end{equation}
		Arguing in components we get
		\begin{align*}
			(\D  w  \cdot Z ) : \Hz  : (\D v \cdot Z) &= (\D  w )_{ip} Z_{pj} \Hz_{ijkl}  (\D v)_{kq} Z_{ql}   =  (\D  w )_{ip}  Z_{pj} \Hz^t_{jikl} Z^\top_{lq} (\D v)_{kq} \\
			&=  (\D  w )_{ip} ((Z \cdot \Hz^t)^t \cdot Z^\top)_{ipkq} (\D v)_{kq} = \D  w  : ( Z\cdot \Hz^t \cdot Z^\top  ) : \D v,
		\end{align*}
		as $(Z \cdot \Hz^t)^t = Z\cdot \Hz^t$, see Lemma~\ref{lem:symmetry of tensor} below. Hence, by setting 
		\begin{equation*}
			\mathbb{C} =  Z\cdot \Hz^t \cdot Z^\top,
		\end{equation*}
		we write
		\begin{equation}\label{eq:derivative of limiting energy in terms of H3}
			\la \partial I( w ) , v \ra =		 \int_\Omega \D  w  :  \mathbb{C} : \D v  \dd x.
		\end{equation}
		Note that major symmetry carries over from $\Hz$ to $\C$.  We also have the following.

		\begin{lemma}[Symmetry of $\C$]\label{lem:symmetry of tensor}
			The tensor $\C$ is minor symmetric, i.e., $\C^t = \C$, and it holds that $(Z \cdot \Hz^t)^t = Z\cdot \Hz^t$. 
		\end{lemma}
		\begin{proof}
			Due to major  symmetry of $\Hz$, it suffices to show that $\C \colon S = 0$ for all skew symmetric  $ S   \in \R^{d \times d}_{\rm skew}$, i.e., $S^\top = -S$. Fix $S \in \R^{d \times d}_{\rm skew}$. Then, $e^{hS} \in \SO(d)$  for $h \in \R$  and by Assumption~\ref{assumptions}   (ii),  we have that 	$0 =  DW(e^{hS}\cdot  Z)$.  Thus, 
			\begin{align*}
				0 &=\frac{\dd}{\dd h} DW (e^{hS} \cdot Z) \vert_{h = 0} = D^2W(e^{hS} \cdot Z) : (e^{hS} \cdot S \cdot Z) \vert_{h = 0} \\
				&= \Hz :(S\cdot Z) =( \Hz \cdot Z^\top ): S.
			\end{align*}
			This   shows $\C \colon S = Z \cdot  \Hz^t \cdot Z^\top  \colon S = 0$. To see the second statement, we use that   $ 0 = (S \cdot Z ) :  \Hz$ by major symmetry of $\Hz$. Therefore, $0 = S : (Z \cdot \Hz^t)^t = S^\top : (Z \cdot \Hz^t)$ for all $S \in \R^{d \times d}_{\rm skew}$ which indeed yields $(Z \cdot \Hz^t)^t = Z\cdot \Hz^t$. 
		\end{proof}

		 Relations   \eqref{eq:derivative of limiting energy in terms of H}--\eqref{eq:derivative of limiting energy in terms of H3}  eventually show 
		\begin{equation}\label{eq: forrrlater}
			\begin{aligned}
				\la \partial I( w ) , v \ra  =  \int_\Omega   (\D  w  \cdot Z) :\mathbb{H}: (\D v \cdot Z) \dd x =  \int_\Omega \D^s  w   :\C : \D^s v \dd x.
			\end{aligned}
		\end{equation}		
		
		\begin{remark}[$L^2$-regularity of $\partial I (w)$]
			If  $\partial I( w )$ is in  $L^2(\Omega;\R^d)$, we have 
			\begin{align}\label{eq:derivative of limiting energy}
				( \partial I( w ) , v ) 	&= -  \int_\Omega \div(\D  w  : \mathbb{C} ) \cdot v \dd x   = -  		\int_\Omega \div(\C :\D^s  w ) \cdot v \dd x, 
			\end{align}	 	
			 	that is, $ - \div (\C: \D^s w) = \partial I(w)$ is the Fr\'echet differential of $I$ with respect to the $L^2$-topology.   Note that this $L^2$-regularity is available  in the purely viscous setting $\rho = 0$, see Remark~\ref{rem:inertia grad in L2} below, or when assuming higher regularity on the initial data, see Lemma~\ref{lem:chain rule} below. 
		\end{remark}

		\subsubsection{ Approximation by scaled  lattice  displacements}\label{subsec:scaled lattice displacements} 
		We start by an auxiliary result  which guarantees the existence of a sequence of scaled lattice displacements converging to a given smooth test function. 
		
		\begin{lemma}[Approximation by scaled lattice displacements]\label{lem:approx by lattice displacements}
			\sloppy	For each    $ v \in C_\comp^{\infty}(\Omega; \R^d)$,  there exists a sequence $(v_k)_k$ of  admissible  scaled lattice displacements  $(v_k)_k$, $  v_k \in   \A_{\epsi_k}$,   such that $v_k \acto v$,  $\Vert v_k \Vert_k \leq \Vert v \Vert_{L^2(\Omega)}$,  and 				
			\begin{equation}\label{eq:L infty bound of disc grad}
				\sup_k \Vert \bar{\D} v_k \Vert_{L^\infty(\ti{\Omega})}  \leq c 	\Vert \D v \Vert_{L^\infty(\Omega)} < \infty,
			\end{equation} 
			 where $c$ depends on $d$ and $A$ only. 
			Moreover,  it holds that  $\bar{\D} v_k \to \D v \cdot Z$ in $L^2(  \Omega; \R^{d \times 2^d})$. 
		\end{lemma}
		
		\begin{proof}
			Extend $v$ to  an element of  $C_\comp^{\infty}(\R^d; \R^d)$ and recall  definition \eqref{eq:definition P} of $P_{\epsi_k}$. By letting the sequence of scaled lattice displacements $v_{k} \colon \L_{\epsi_k} \cap \ti{\Omega} \to \R^d$ be defined by 
			\begin{align*}
				v_{k} (x) := \left\lbrace   \begin{array}{lcl}
					P_{\epsi_k} v(x) & \text{for  }   &  x \in \L_{\eps_k} \cap \Omega,  \\
					0 & \text{for  }   &  x \in \L_{\eps_k} \cap (\ti{\Omega}\setminus \Omega),
				\end{array} \right.
			\end{align*}
			the convergence $v_k \acto v$ follows directly by definition \eqref{eq:convergence}.   Moreover, $\Vert v_k \Vert_k \leq \Vert v \Vert_{L^2(\Omega)}$  is a consequence of  the definition of $v_k$,  \eqref{eq:atomistic inner product},  and Jensen's inequality.

 We now show \eqref{eq:L infty bound of disc grad}. 			We fix $\bar{x} \in \L_{\epsi_k}' \cap \ti{\Omega}$ and write  $v_k^i = v_k^i(\bar{x})   = v_k (\ove{x} + \epsi_k z_i)$, $i = 1, \dots, 2^d$,  and $\bar{v}_k    = \bar{v}_k(\bar{x}) = 2^{-d} \sum_{i=1}^{2^d}  v_k^i$.  Let $j_0 \in \{1, \dots, 2^d\}$ be such that $\max_{j = 1, \dots, 2^d} \vert v_k^i(x) - v_k^j (x)\vert = \vert v_k^i(x) - v_k^{j_0} (x) \vert$.  To obtain an $L^\infty$-bound on the discrete gradient, we use the Mean-Value Theorem to compute  
			\begin{align*}
				\vert v_{k}^i - \bar{v}_k \vert \leq  \vert v_k^i - v_k^{j_0} \vert \leq c\epsi_k \Vert \D v \Vert_{L^\infty(\ti{\Omega})}  \leq c \epsi_k \Vert \D v \Vert_{L^\infty(\Omega)},  
			\end{align*} 
			where $c = c(d,A)$ is the diameter of  $ \epsi_k^{-1}  A[\bar{x} -\frac{3}{2}\epsi_k, \bar{x} +\frac{3}{2} \epsi_k]^d$, noting that $Q_\eps(\ove{x} + \epsi_k z_i) \subset A[\bar{x} -\frac{3}{2}\epsi_k, \bar{x} +\frac{3}{2} \epsi_k]^d$ for all $i = 1, \dots, 2^d$.  			As the discrete gradient is constant on individual scaled cells, we obtain
			\begin{align*}  
				\Vert \bar{\D} v_k \Vert_{L^\infty(\ti{\Omega})} &= 			\sup_{ \bar{x} \in  \L'_{\eps_k}(\ti{\Omega})} \Big( \sum_{ i  = 1}^{2^d}  \frac{1}{\epsi_k^2}  \vert v_k^i ( \bar{x}  ) - \bar{v}_k( \bar{x}) \vert^2 \Big)^{1/2}  \leq \big(  2^d  c^2 \Vert \D v \Vert_{L^\infty(\Omega)}^2 \big)^{1/2} < \infty. 
			\end{align*} 
			For the proof of $\bar{\D} v_k \to \D v \cdot Z$ in $L^2( \Omega;\R^{d \times 2^d})$, we  refer  to \cite[Lemma~4.4]{Schmidt_atomisticToCont}. 
		\end{proof}

		Let us now introduce the \df{piecewise constant interpolation} $\ti{u}_k$ on scaled cells of an admissible scaled lattice displacement $u_k \in \A_{\epsi_k}$.   To this end, recall the notation  $\bigcup Q_{\epsi_k} = \bigcup_{\ove{x} \in \L'_{\epsi_k}(\ti{\Omega})} Q_{\epsi_k}(\ove{x})$. For $\xi \in \bigcup Q_{\epsi_k}$, let  $x \in \L_{\epsi_k} \cap \ti{\Omega}$ be such that $\xi \in Q_{\epsi_k}( x  )$. Then, the piecewise constant interpolation is defined as 
		\begin{equation}\label{eq:pw const int}
			\ti{u}_k(\xi) := \left\lbrace \begin{array}{lcl}
				u_k(x)	& \text{if} & Q_{\epsi_k}( x  ) \cap \Omega\neq \emptyset \\
				0 & \text{if}  & Q_{\epsi_k}( x  ) \cap \Omega = \emptyset.
			\end{array}\right.
		\end{equation}
		Moreover, we trivially extend $\ti{u}_k$ to the whole of $\ti{\Omega}$  without changing notation.  We thus  get  for $\epsi_k > 0$ small enough that
		\begin{align}\label{eq:sum and int}
			\epsi_k^d \det A  \sum_{ x \in  \L_{\epsi_k} \cap \ti{\Omega}} u_k(x) &= 	\epsi_k^d \det A  \sum_{ x \in  \L_{\epsi_k} \cap \ti{\Omega}} \fint_{Q_{\epsi_k}(x) } \ti{u}_k(\xi) \dd\xi \nonumber \\
			&=  \int_{\bigcup Q_{\epsi_k}}  \ti{u}_k (\xi ) \dd \xi = \int_{\ti{\Omega}} \ti{u}_k (\xi) \dd \xi.
		\end{align} 
		In particular,  by a similar argument  we  also  have that 
		\begin{align}\label{eq:norms are the same}
			\Vert u_k \Vert_k = \Vert \ti{u}_k \Vert_{L^2(\ti{\Omega})}.		
		\end{align} 
				The above interpolation is motivated by the following properties.

		\begin{lemma}[Equivalence of atomistic and  $L^2$-convergence]\label{lem:conv of norms}
			 Letting  $(u_k)_k$, $ u_k \in   \mathcal{A}_{\epsi_k} $, be a sequence of  admissible  scaled lattice displacements and $u \in L^2(\Omega;\R^d)$,   then 
			\begin{align*}
				u_k \acto u \quad \iff \quad \int_{\ti{\Omega}} \vert \ti{u}_k - \ti{u} \vert^2 \to 0,
			\end{align*} 
			where $\ti{u}$ is the trivial extension of $u$ to $\ti{\Omega}$. 
			Moreover, for $u,v \in L^2(\Omega;\R^d)$ and $(u_k)_k$, $(v_k)_k$ with $u_k, v_k \in \A_{\epsi_k}$ such that $u_k \acto u$ and $v_k \acto v$, it holds that
			\begin{align*}
				(u_k, v_k)_k \to (u,v).
			\end{align*}
			In particular, we have that $\Vert u_k \Vert_k  \to \Vert u \Vert_{L^2(\Omega)}$ for $k \to \infty$.
		\end{lemma}

		\begin{proof}
			At first, note that AC-convergence implies strong $L^2$-convergence of piecewise constant interpolations. Indeed, let $\ti{u}_k$ and $\ti{P}_{\epsi_k} u$ be the piecewise constant interpolations of $u_k$ and $P_{\epsi_k}u$, respectively. Then, extending $u$ trivially as $\ti{u}$ to $\ti{\Omega}$, we obtain
			\begin{align*}
				\int_{\ti{\Omega}} | \ti{u}_k (\xi ) - \ti{u}(\xi) |^2 \dd \xi  \leq 2 \int_{\ti{\Omega}} | \ti{u}_k (\xi) - \ti{P}_{\epsi_k}u(\xi) |^2 \dd \xi + 2 \int_{\ti{\Omega}} | \ti{P}_{\epsi_k}u(\xi) - \ti{u}(\xi) |^2 \dd \xi \\
				= 2 \epsi_k^d \det A  \sum_{ x \in  \L_{\epsi_k} \cap \ti{\Omega}} \vert u_k(x) - P_{\epsi_k}u(x) \vert^2 + 2\int_{\ti{\Omega}} | \ti{P}_{\epsi_k}u(\xi) - \ti{u}(\xi) |^2 \dd \xi \to 0.
			\end{align*}
			Here, we used that $\ti{P}_{\epsi_k}u\to \ti{u}$ strongly in $L^2(\ti{\Omega};\R^d)$.    Conversely, strong $L^2$-convergence of piecewise constant interpolations implies AC-convergence   as 
			\begin{align*}
				\epsi_k^d \det A  \sum_{ x \in  \L_{\epsi_k} \cap \ti{\Omega}} \vert u_k(x) - P_{\epsi_k}u(x) \vert^2   &= \int_{\ti{\Omega}} \vert \ti{u}_k - \ti{P}_{\epsi_k}u \vert^2 \\
				&\leq 2\int_{\ti{\Omega}} \vert \ti{u}_k - \ti{u} \vert^2 +2  \int_{\ti{\Omega}} \vert  \ti{u}  - \ti{P}_{\epsi_k}u \vert^2 \to 0.
			\end{align*}
			Eventually, due to $\ti{u}_k \to u$ and $\ti{v}_k \to v$ in $L^2(\ti{\Omega};\R^d)$,  we get that 
			\begin{align*}
				(u_k , v_k )_k &= 	\epsi_k^d \det A \sum_{ x \in  \L_{\epsi_k} \cap \ti{\Omega}} u_k(x) \cdot v_k(x) \\
				&= \int_{ \ti{\Omega}} \ti{u}_k (\xi) \cdot \ti{v}_k (\xi) \dd \xi \to \int_\Omega u(\xi) \cdot v(\xi) \dd \xi = (u,v). \qedhere
			\end{align*}
		\end{proof}
		 
                \subsection{A priori estimates}\label{sec:apriori}
Let  $u_k$ be the solution of the atomistic equation of motion
\eqref{eq:atomistic the second}.  By using the chain rule,  for all $t\in [0,T]$ we get
\begin{align*}
  0 &  = \frac{1}{2\nu} \int_0^t \| \rho \ddot u_k(s) + \nu \dot u_k(s) +
      \partial I_k(u_k(s))\|_k^2 \, \dd s \\
  &=\int_0^t \left( \big( \rho \ddot{u}_k(s) + \partial I_k (u_k(s)) ,
    \dot{u}_k(s) \big)_k + \frac{\nu}{2}\|\dot u_k(s)\|^2_k +
    \frac{1}{2\nu} \| \rho \ddot u_k (s) +
    \partial I_k(u_k(s))\|_k^2 \right) \, \dd s\\
   &= \frac{\rho}{2}\|\dot u_k(t)\|^2_k + I_k(u_k(t)) -
     \frac{\rho}{2}\| \dot u_k(0)  \|^2_k - I_k( u_k(0) ) \\
  &\quad + \frac{\nu}{2}\int_0^t   \|\dot u_k(s)\|^2_k \, \dd s +
    \frac{1}{2\nu} \int_0^t \| \rho \ddot u_k (s) +
    \partial I_k(u_k(s))\|_k^2  \, \dd s,
\end{align*}
  where   $\partial I_k (u_k(t))$ denotes the scaled lattice  strain given in \eqref{eq:defin disc grad as lattice disp}.
In particular, given the initial conditions $u_k(0)=u^0_k  := u^0_{\epsi_k}$ and $\dot u_k(0)=u^1_k   := u^1_{\epsi_k}$, system \eqref{eq:atomistic the second} is equivalent to
the energy-dissipation-inertia equation
\begin{align}\label{eq:EDIE}
  &\frac{\rho}{2}\|\dot u_k(t)\|^2_k + I_k(u_k(t)) + \frac{\nu}{2}\int_0^t   \|\dot u_k(s)\|^2_k \, \dd s +
    \frac{1}{2\nu} \int_0^t \| \rho \ddot u_k (s) +
    \partial I_k(u_k(s))\|_k^2  \, \dd s \nonumber \\
  &\quad =  \frac{\rho}{2}\|u^1_k\|^2_k + I_k(u^0_k) \quad
    \forall t\in [0,T],
\end{align}
or, equivalently, for $t=T$.  
As we are assuming  that  $u^1_k \acto w^1 \in L^2(\Omega;\Rz^d)$ and  that 
$I_k(u^0_k)$ is bounded, from the energy-dissipation-inertia equation
we obtain the following.

\begin{lemma}[A priori estimates]\label{lem:apriori}
	Let $u_k$ be the solution of the atomistic equation of motion
  \eqref{eq:atomistic the second}. Then,
  \begin{align}
    \label{eq:apriori}
    \sqrt{\rho}   \max_{t\in[0,T]}\| \dot u_k(t)\|_k&\leq c,\\                                          
      \max_{t\in[0,T]}I_k(u_k(t)) &\leq c, \label{eq:uniform energy}\\
     \int_0^T\| \dot u_k(s)\|_k^2 \, \dd s &\leq c,  \label{eq:H1L2 bound 1}\\
    \int_0^T \| \rho \ddot u_k(s) + \partial I_k(u_k(s))\|_k^2 \, \dd s &\leq c, \label{eq:unform bound on force}
  \end{align}
  where $c > 0 $ depends  only  on  $  \rho  \sup_k \|u^1_{ k}\|^2_{ k}$,
  $\sup_{ k}   I_k( u^0_k  ))$,   and $\nu$.
\end{lemma}

		\subsection{Compactness and limit of $\nu \dot{u}_k$}\label{subsec:limit2}
				
		 	In this section, we show that the sequence $(u_k)_k$ of solutions to the atomistic problem~\eqref{eq:atomistic the second} admits an AC-converging subsequence $u_k(t) \acto w(t) \in L^2(\Omega;\R^d)$ for all $t \in [0,T]$. Moreover, as a first step towards showing that the limit $w$ solves the continuum problem, we also prove that $\nu \dot{u}_k \wto \nu \partial_t w$ in a weak $L^2(0,T;L^2(\Omega;\R^d))$-sense made precise below.

	 From  \eqref{eq:H1L2 bound 1} we obtain that  $\dot{u}_k$ is  uniformly  bounded in $L^2(0,T; (\A_{\epsi_k}, \Vert \cdot \Vert_k))$. Let now  $\partial_t \ti{u}_k$  indicate the piecewise constant interpolants of $\dot{u}_k$  from \eqref{eq:pw const int}.  In view of \eqref{eq:norms are the same},   we get  
		\begin{align}\label{eq:H1L2 bound 2}
			\int_0^{ T} \Vert \dot{u}_k (s, \cdot )\Vert_k^2 \dd s = \int_0^{ T} \Vert \partial_t \ti{u}_k (s , \cdot) \Vert^2_{L^2(\ti{\Omega}) } \dd s  \le c. 
		\end{align}		
 This along with the compactness result in  Theorem~\ref{thm:gamma convergence and compactness} yields that the solution $u_k$ to the atomistic problem \eqref{eq:atomistic the second} satisfies that $(\ti{u}_k(t))_k$ is compact in $L^2(\ti{\Omega};\R^d)$ for all $t \in (0,T)$. This leads to the  following.

		\begin{lemma}[Compactness]\label{lem:Aubin Lions}
				 Let $(u_k)_k$ be a sequence of solutions  to the atomistic problem \eqref{eq:atomistic the second}. Then, there exists  $w \in C([0,T], L^2(\Omega; \R^d))  \cap L^\infty(0,T;H^1_0(\Omega;\R^d))$ and a (non relabeled) subsequence of $(u_k)_k$ such that
				\begin{align*}
					\ti{u}_k \to \ti{w} \quad \text{in  $C([0,T], L^2(\ti{\Omega};\R^d))$,} 
				\end{align*}
				where $\ti{u}_k$ is the piecewise constant interpolation of $u_k$ and $\ti{w}$ denotes the trivial extension of $w$ to $\ti{\Omega}$. 
		\end{lemma}
		\begin{proof}
The proof relies on a variant of the Aubin-Lions Theorem. However, due to the discrete nature of our spaces, it is more convenient to directly resort to the  Arzel\`a-Ascoli  Theorem. 
			First, note that due to \eqref{eq:uniform energy} and Theorem~\ref{thm:gamma convergence and compactness} the set $\{ \ti{u}_k (t) : k \in \N \}$ is relatively compact in $L^2(\ti{\Omega};\R^d)$ for each $t \in [0,T]$. Equi-continuity of $(\ti{u}_k)_k$ follows from   \eqref{eq:H1L2 bound 2}. Indeed,    \eqref{eq:H1L2 bound 2} implies   that $(\ti{u}_k)_k \subset H^1(0,T;L^2(\ti{\Omega};\R^d)) \subset C([0,T], L^2(\ti{\Omega};\R^d))$ and for each $ 0 < s\le t <T$
			\begin{align*}
				\Vert \ti{u}_k(t) - \ti{u}_k(s) \Vert_{L^2(\ti{\Omega})} \leq \int_s^t \Vert \partial_t \ti{u}_k (\tau) \Vert_{L^2(\ti{\Omega})} \dd \tau \leq  \sqrt{t-s} \Vert \partial_t \ti{u}_k \Vert_{L^2(0,T;L^2(\ti\Omega;\R^d))} \le c\sqrt{t-s}. 			
			\end{align*}
		  	Thus,   $(\ti{u}_k)_k$  is equi-continuous and    the convergence to a limit  $w \in C([0,T], L^2(\ti{\Omega};\R^d))$  follows   from the Arzel\`a-Ascoli Theorem.  Since $u_k$ is extended by $0$ on $\L_{\epsi_k} \cap (\ti{\Omega} \setminus \Omega)$, the limit $w$ is indeed supported on $\Omega$.  Eventually, applying again \eqref{eq:uniform energy} and Theorem~\ref{thm:gamma convergence and compactness}  we find that $w(t) \in H^1_0(\Omega)$ for all $t \in [0,T]$ with uniformly controlled $H^1_0(\Omega)$-norm. 
		\end{proof}

		The lemma above  implies   that $u_k(t) \acto w(t)$ for all $ t \in [0,T]$ as $L^2$- and AC-convergence are equivalent due to Lemma~\ref{lem:conv of norms}. 
		
		 We now  prove that $w$ is indeed a solution to the continuum problem \eqref{eq:continuous problem}. 		As a first step, we show that   the piecewise constant interpolations of   $\nu \dot{u}_k$   converge weakly to  $\partial_t w $  in $L^2(0,T; L^2(\ti{\Omega};\R^d))$. 	Indeed,   from \eqref{eq:H1L2 bound 2}  we can  extract a (not
		relabeled) subsequence   of $u_k$ such that  		 $\partial_t \ti{u}_k \wto f$ in $L^2(0,T;L^2(\ti{\Omega};\R^d))$.  From the linearity of the time derivative 	we deduce that $f = \partial_t \ti{ w }$, where $\ti{ w }$  again denotes  the trivial extension of $ w $ to $\ti{\Omega}$.   Thus, we have shown  
		\begin{align}\label{eq:weak conv dot u}
			\partial_t \ti{u}_k \wto    \partial_t \ti{ w }     \quad \text{ 	 in
				$L^2(0,T;L^2(\ti{\Omega};\R^d))$. }
		\end{align}
 	
			\begin{remark}[Purely viscous case $\rho = 0$] \label{rem:inertia grad in L2-again}  			
				 Note that the results of this section stay valid when considering the purely viscous case of $\rho = 0$. Indeed, as also the energy-dissipation equality \eqref{eq:EDIE} holds (along with the obvious modifications), Lemmas~\ref{lem:apriori} and \ref{lem:Aubin Lions} remain true,  where \eqref{eq:apriori} becomes trivial.  
\end{remark}

		\subsection{Limit of $\rho \ddot{u}_k + \partial I_k(u_k)$}\label{subsec:limit1}
		 In this section, we prove that,  given   the pointwise in time AC-converging sequence $u_k \acto w$  from Lemma~\ref{lem:Aubin Lions},  we also have that the piecewise constant interpolants of    $\rho \ddot{u}_k + \partial I_k(u_k)$  converge weakly to     $\rho \partial_{tt} w + \partial I(w)$   in $L^2(0,T;L^2( {\Omega};  \R^d))$.

	  
		Denoting by $	(\rho \ddot{u}_k + \partial I_k(u_k))^\sim $ the  interpolation of $ \rho \ddot{u}_k + \partial I_k(u_k)$  introduced in \eqref{eq:pw const int},   \eqref{eq:unform bound on force} together with \eqref{eq:norms are the same}     implies that there exists $\alpha \in L^2(0,T;L^2( \ti{\Omega};  \R^d))$ such that   
		\begin{equation}\label{eq:conv to alpha}
				(\rho \ddot{u}_k + \partial I_k(u_k))^\sim \wto \alpha \quad \text{  weakly   in  } L^2(0,T;L^2(\ti{\Omega}; \R^d)).
		\end{equation}
		  We are left with showing that $\alpha = \rho \partial_{tt}  \ti{w}  + \partial I(\ti{w})$ in $L^2(0,T;L^{2}( \ti{\Omega};  \R^d))$,  where $\ti{w}$ is the trivial extension of $w$.  
		 For the identification of the limit, we  separately  investigate the limits of  $\rho \partial_{tt}\ti{u}_k(t)$ and $\left(\partial I_k(u_k(t))\right)^\sim$, i.e., the limits of the piecewise constant interpolations of $\rho\ddot{u}_k$ and $\partial I_k(u_k)$.  This is a delicate issue, for  $\partial I_k (u_k)$ is defined on a space  which depends on $k$.  As an important ingredient, we    deduce separate   uniform bounds on $\rho \partial_{tt}\ti{u}_k(t)$ and $\left(\partial I_k(u_k(t))\right)^\sim$, respectively,  by resorting to a sufficiently weak topology.  	   We treat the respective limits  separately  in Sections \ref{sec: XXX1}--\ref{sec: XXX2}, and identify $\alpha$ in Section \ref{sec: XXX3}.

		\subsubsection{Limit of $\partial I_k (u_k)$}\label{sec: XXX1}  Here, we prove convergence of    $\partial I_k (u_k)$ to $\partial I(w)$ in the sense of  G\^{a}teaux derivatives.  The analysis  of  this  section  is carried out for $t \in (0,T)$ arbitrary but fixed.  Correspondingly, we omit to indicate the  dependence on $t$ for  brevity.  We start by computing the variation of $I_k$ at $u_k$ in direction $v_k$ as
		\begin{align*}
			\dfrac{\dd}{\dd h} I_k (u_k + h v_k ) \big\vert_{h = 0} = 	\eps_k^d \det A \sum_{\overline{x} \in \L'_{\epsi_k}(\ti{\Omega})} \frac{1}{\delta_k} D W (Z + \delta_k \bar{\D} u_k(\ove{x})) \cdot \bar{\D} v_k(\ove{x}).
		\end{align*} 
		We  have the following.

			\begin{lemma}[Convergence of G\^{a}teaux derivatives]\label{lem:convergence} 
			Let $(u_k)_k$ be a sequence of admissible scaled lattice displacements, $u_k \in  \mathcal{A}_{\epsi_k}$, with equibounded atomistic energy, i.e., $\sup_k I_k(u_k) < \infty$.  Suppose  that  $u_k \acto  w $  with $w \in H^1_0(\Omega;\R^d)$.   	Assuming $v \in \Cc^\infty_\comp(\Omega; \R^d)$ with $v_k$ being the corresponding scaled lattice displacements from Lemma~\emph{\ref{lem:approx by lattice displacements}}, we have 
			\begin{align*}
				\epsi_k^d  \det A  \sum_{\ove{x} \in  \L'_{\epsi_k} (\ti{\Omega})} \frac{1}{\delta_k} DW(Z + \delta_k \bar{\D} u_k (\ove{x})) : \bar{\D} v_k(\ove{x})	 \to \int_\Omega (\D  w  \cdot Z ) : \Hz : (\D v \cdot Z) \dd x
			\end{align*}
			as $k \to \infty$, where we recall that the right-hand side coincides with   $\la \partial I(w), v\ra$, see \eqref{eq:derivative of limiting energy in terms of H}. 
		\end{lemma}
		\begin{proof}
			We compute the Taylor expansion of $D W$ at $Z$. This calls for $\delta_k \bar{\D} u_k$ to be small. We hence define  \df{good sets} $G_k$ and  \df{bad sets} $B_k$ as
			\begin{align*}
				G_k &:= \left\{ \ove{x} \in \L'_{\epsi_k}(\ti{\Omega}) :  \vert \delta_k \bar{\D} u_k (\ove{x}) \vert < \beta_k^{-1}   \right\},  \\
				B_k &:= \left\{ \ove{x} \in \L'_{\epsi_k}(\ti{\Omega}) :  \vert \delta_k \bar{\D} u_k (\ove{x}) \vert \geq \beta_k^{-1}   \right\},
			\end{align*}
			where $\beta_k$ is chosen in such a way that $\beta_k \to \infty$, and $\beta_k^2 \delta_k \to 0$ as $k \to \infty$.
 			Denoting by $\# B_k$ the cardinality of $B_k$, a Chebyshev argument ensures that   $\delta_k^{-1}  \eps_k^d  \# B_k \to 0$  for $k\to \infty$. Indeed,  by   relation \eqref{eq:norms are the same} we estimate 
			\begin{equation}\label{eq:B_k vanishes}
				\frac{ \eps_k^d }{\delta_k} \# B_k = \frac{ \eps_k^d  \beta_k^2}{\delta_k} \sum_{\ove{x} \in B_k} \frac{1}{\beta_k^2} \leq   \frac{ \eps_k^d \beta_k^2}{\delta_k} \sum_{\ove{x} \in B_k}  \delta_k^2 \vert \bar{\D} u_k(\ove{x}) \vert^2 \leq \frac{\beta_k^2 \delta_k}{ {\rm det} A } \Vert \bar{\D} u_k \Vert_{L^2(\ti{\Omega})}^2 \to 0,
			\end{equation}
			 where we have used $\sup_k \Vert \bar{\D} u_k \Vert_{L^2(\ti{\Omega})}  < \infty$,   which  follows from the compactness result of Theorem~\ref{thm:gamma convergence and compactness}. A fortiori, we have $  \eps_k^d  \# B_k \to 0$ as $k \to \infty$.
			
			Now, we Taylor-expand $DW$ at $Z$ for $\ove{x} \in G_k$ obtaining
			\begin{align*}
				D W (Z + \delta_k \bar{\D} u_k (\ove{x}) ) &= D^2 W (Z) : (\delta_k \bar{\D} u_k (\ove{x}) ) + R( \delta_k \bar{\D} u_k (\ove{x}) ) \\
				&=  \delta_k \mathbb{H} : \bar{\D} u_k (\ove{x}) +   R( \delta_k \bar{\D} u_k (\ove{x}) ),
			\end{align*}
			where the  \df{remainder term}  $R \colon \R^{d\times 2^d} \to \R^{d\times 2^d}$  satisfies \[ \sup \left\{  \frac{  \vert R( F ) \vert  }{\vert F \vert } \, : \, \vert  F \vert \leq \rho  \right\} \to 0 \quad \text{   as  } \rho \to 0. \]  Considering only the energy of cells in $G_k$,  due to  the  major  symmetry of $\Hz$,  we arrive at
			\begin{align}\label{eq:int and sum}
				&\epsi_k^d \det A \sum_{\overline{x} \in \L'_{\epsi_k}(\ti{\Omega}) \cap G_k} \frac{1}{\delta_k} D W (Z + \delta_k \bar{\D} u_k(\ove{x})) : \bar{\D} v_k (\ove{x}) \notag \\
				&= \epsi_k^d \det A \sum_{\overline{x} \in \L'_{\epsi_k}(\ti{\Omega}) \cap G_k} \left(\bar{\D} u_k (\ove{x}) :  \Hz :  \bar{\D} v_k (\ove{x})  + \delta_k^{-1} \bar{\D} v_k  (\ove{x}) :   R( \delta_k \bar{\D} u_k(\ove{x}) ) \right).
			\end{align}
			The next step is to control $\delta_k^{-1} \bar{\D} v_k : R( \delta_k \bar{\D} u_k)$ on $G_k$.   Note that, as $\bar{\D} u_k(\ove{x})$ is constant on scaled cells $Q_{\epsi_k}(\ove{x})$, we have,  analogously to \eqref{eq:sum and int},  $ \int_{Q_{\epsi_k}(\ove{x})} \bar{\D} u_k (x) \dd x = \epsi_k^d \det A \, \bar{\D}u_k(\ove{x})$, and thus, 
			\begin{equation}\label{eq: LLL}
				\int_{\bigcup Q_{\epsi_k}} \bar{\D} u_k (x) \dd x = \epsi_k^d \det A \sum_{\ove{x} \in \L'_{\epsi_k}(\ti{\Omega})} \bar{\D} u_k(\ove{x}),
			\end{equation}
			where we 				recall    that   $\bigcup Q_{\epsi_k} = \bigcup_{\ove{x} \in \L'_{\epsi_k}(\ti{\Omega})} Q_{\epsi_k}(\ove{x})$. 
			Therefore, letting $k \to \infty$,  by H\"older's inequality  we get
			\begin{align}\label{eq:little o vanishes}
				&\epsi_k^d  \det A   \Bigg|  \sum_{\overline{x} \in \L'_{\epsi_k}(\ti{\Omega}) \cap G_k}  \delta_k^{-1} \bar{\D} v_k (\ove{x}): R( \delta_k \bar{\D} u_k(\ove{x})) \Bigg|   \nonumber \\
				&\leq  \int_{ \bigcup_{G_k} Q_{\epsi_k} } \vert  \bar{\D} v_k(x) \vert \  \vert \bar{\D} u_k(x) \vert  \frac{  \vert R(\delta_k \bar{\D} u_k (x))  \vert   }{\vert \delta_k \bar{\D} u_k (x) \vert }  \dd x \nonumber \\ 
				&\leq  \sup_{\bigcup_{G_k} Q_{\epsi_k}} \left(\frac{  \vert  R(\delta_k\bar{\D} u_k )  \vert  }{\delta_k \vert \bar{\D} u_k \vert}\right)    \Vert \bar{\D} v_k \Vert_{L^2(\ti{\Omega})}   \Vert\bar{\D} u_k \Vert_{L^2(\ti{\Omega})}  \to 0,
			\end{align}
			where we set $\bigcup_{G_k} Q_{\epsi_k} := \bigcup_{\ove{x} \in G_k \cap \L'_{\epsi_k}(\ti{\Omega})} Q_{\epsi_k}(\ove{x})$.			
			Convergence \eqref{eq:little o vanishes} follows as $\vert \delta_k \bar{\D} u_k \vert < \beta_k^{-1}$  uniformly  on $G_k$ with $\beta_k^{-1} \to 0$. Hence, the  above supremum  over $\bigcup_{G_k} Q_{\epsi_k}$ converges to $0$.  Recall  also  that $\Vert \bar{\D} v_k \Vert_{ L^2(\ti{\Omega})}$ and $ \Vert\bar{\D} u_k \Vert_{L^2(\ti{\Omega})} $ are uniformly bounded by Lemma~\ref{lem:approx by lattice displacements} and Theorem~\ref{thm:gamma convergence and compactness}, respectively. 
			
			On the  bad sets  $B_k$, we use Assumption~\ref{assumptions} (v) and the  fact that  $\Vert \bar{\D} v_k \Vert_{L^\infty(\ti{\Omega})}  \le c$, see  Lemma~\ref{lem:approx by lattice displacements}, to obtain
		\begin{align}\label{eq:sum on B_k vanishes}
			&\frac{\epsi_k^d \det A}{\delta_k} \sum_{\overline{x} \in \L'_{\epsi_k}(\ti{\Omega}) \cap B_k } D W (Z+  \delta_k \bar{\D} u_k) :  \bar{\D} v_k \leq  \frac{1}{\delta_k} \int_{\bigcup_{B_k} Q_{\epsi_k}} \vert D W (Z + \delta_k \bar{\D} u_k)\vert \, \vert \bar{\D} v_k \vert \dd x \nonumber \\
			&\leq \frac{ c  }{\delta_k} \Vert \bar{\D} v_k \Vert_{L^\infty(\ti{\Omega})} \int_{\bigcup_{B_k} Q_{\epsi_k}} \vert Z + \delta_k \bar{\D} u_k \vert^2  \dd x + \frac{ c  }{\delta_k}  {\rm det} A \, \eps_k^d  \Vert \bar{\D} v_k \Vert_{L^\infty(\ti{\Omega})}  \# B_k  \nonumber \\
			&\leq \frac{c}{\delta_k} \int_{\bigcup_{B_k} Q_{\epsi_k}} \left( |Z|^2+ \delta_k^2 \vert \bar{\D} u_k \vert^2 \right) \dd x  + \frac{ c }{\delta_k}  {\rm det} A \, \eps_k^d   \# B_k  \nonumber \\
			& \leq \frac{c}{\delta_k}  {\rm det} A \,  \eps_k^d  \# B_k  |Z|^2   + c \delta_k\Vert \bar{\D} u_k \Vert_{L^2(\ti{\Omega})}^2  + \frac{ c }{\delta_k}  {\rm det} A \, \eps_k^d   \# B_k  \to 0,
		\end{align}
		where we set $\bigcup_{B_k} Q_{\epsi_k} := \bigcup_{\ove{x} \in B_k \cap \L'_{\epsi_k}(\ti{\Omega})} Q_{\epsi_k}(\ove{x})$ and used the fact that  $\delta_k^{-1}  \eps_k^d  \# B_k \to 0$  due to \eqref{eq:B_k vanishes}. Here, we again employed Theorem~\ref{thm:gamma convergence and compactness}  to bound $\sup_k \Vert \bar{\D} u_k \Vert_{L^2(\ti{\Omega})}^2 \le c$. 
			
			By Lemma~\ref{lem:approx by lattice displacements} we have that   $\bar{\D} v_k \to \D v \cdot Z$ strongly in $L^2(\ti{\Omega}; \R^{d \times 2^d})$.   Indicating by  $\chi_{_{\bigcup_{G_k} Q_{\epsi_k}} }$  the characteristic function of the set ${\bigcup_{G_k} Q_{\epsi_k}}$, we observe that $\chi_{_{\bigcup_{G_k} Q_{\epsi_k}} }\to 1$ in measure. Recall that $\bar{\D} u_k \wto \D   \ti{w}   \cdot Z$  weakly   in $L^2( \ti{\Omega};\R^{d\times 2^d})$ by Theorem~\ref{thm:gamma convergence and compactness}. By the linearity of $F \mapsto \Hz : F$,  we  obtain 
			\begin{equation*}
				\chi_{_{\bigcup_{G_k} Q_{\epsi_k}} } \:	\Hz: (\bar{\D} u_k) \wto  \Hz : (\D  w  \cdot Z) \text{   weakly   in  } L^2( \ti{\Omega};  \R^{d\times 2^d}). 
			\end{equation*}
			This, together  with \eqref{eq: LLL} and  \eqref{eq:little o vanishes},  implies 
			\begin{align*}
				&\epsi_k^d  \det A \sum_{\overline{x} \in \L'_{\epsi_k}(\ti{\Omega}) \cap G_k} \left(  \bar{\D} u_k :  \Hz :  \bar{\D} v_k  + \delta_k^{-1} \bar{\D} v_k  :   R( \delta_k \bar{\D} u_k) \right) \\
				&\to \int_\Omega (\D  w  \cdot Z) : \Hz : (\D v \cdot Z) \dd x
			\end{align*}
			which, along with  \eqref{eq:int and sum}  and \eqref{eq:sum on B_k vanishes},	 concludes the proof.
		\end{proof}

		\subsubsection{Limit of  $\rho \ddot{u}_k$}\label{sec: XXX2}

		We now show that  $\rho \ddot{u}_k$  converges suitably  to $\rho \partial_{tt} w$. To this end, we uniformly bound  $\rho \ddot{u}_k$  in a weaker topology to obtain compactness and then identify the corresponding limit. 	
		
		Let $\ell > \frac{d}{2} + 1$,  so that   $ H^{\ell}_0(\Omega;\R^d) \subset W^{1,\infty}(\Omega;\R^d)$, see, e.g., \cite[Theorem~4.12]{afSob}. We define  $\ove{\ddot{u}_k} \in L^2(0,T;H^{-\ell}(\Omega;\R^d))$ almost everywhere in time as 
		\begin{align}\label{eq:ove ddot u}
			\la \ove{\ddot{u}_k(t)} , v \ra_{H^{\ell}_0(\Omega)} := (\ddot{u}_k(t) , v_k)_k   \quad \quad \text{for all $v\in H^\ell_0(\Omega;\R^d)$, }  
		\end{align} 
		where $v_k$ are the scaled lattice displacements of $v$ from Lemma~\ref{lem:approx by lattice displacements}. Note that, strictly speaking, we define \eqref{eq:ove ddot u} for $v \in C^\infty_c(\Omega;\R^d)$ only and argue by density, as the approximation by scaled lattice displacements from Lemma~\ref{lem:approx by lattice displacements}  has only been introduced  for smooth functions with compact support.

		Analogously to $\ove{\ddot{u}_k}$ in  \eqref{eq:ove ddot u}, we define  $\overline{ \partial I_k (u_k) } \in L^2(0,T;H^{-\ell}(\Omega;\R^d))$ by setting, again almost everywhere in time, 
		\begin{align}\label{eq:lip op}
			\langle	\overline{ \partial I_k (u_k(t)) } , v \rangle_{ H^{\ell}_0(\Omega) } := \langle \partial I_k (u_k(t)), v_k \rangle    \quad \quad \text{for all $v\in H^\ell_0(\Omega;\R^d)$, }  
		\end{align}
		where the right-hand side of \eqref{eq:lip op} is defined by \eqref{eq:defin of gradient}. We  now  prove that $\ove{\ddot{u}_k(t)}$ and $\ove{\partial I_k(u_k(t))}$   have $L^2$-regularity  in time  uniformly in $k$.

		\begin{lem}[Uniform bound]\label{lem:uniform bound on operator}
	 The  functions      $\rho  \ove{\ddot{u}_k} $ and      $\overline{ \partial I_k (u_k) }$ are uniformly bounded in $L^2(0,T;H^{-\ell}(\Omega;\R^d))$, i.e., 
			\begin{align*}
	{\rm (i)} \ \ 			\sup_k \int_0^T \left\Vert \ove{\partial I_k (u_k(t) )}  \right\Vert_{H^{-\ell}(\Omega)}^2 \dd t \leq c, \quad  \quad { {\rm (ii)} \ \  \sup_k \int_0^T \left\Vert \rho  \ove{\ddot{u}_k(t)} \right\Vert_{H^{-\ell}(\Omega)}^2 \dd t \leq c.   }
			\end{align*}
		\end{lem}

		\begin{proof} 
		 We start by checking (i).  Let $t \in (0,T)$ be fixed and recall  by \eqref{eq:defin of gradient} that 
			\begin{align*}
				\la \ove{\partial I_k(u_k(t))}, v \ra_{ H^{\ell}_0(\Omega)} &= \la \partial I_k (u_k(t, \cdot)) , v_k \ra \\
				&= \epsi_k^d \det A \sum_{\ove{x} \in  \L'_{\epsi} (\ti{\Omega})} \frac{1}{\delta_k} DW(Z + \delta_k \bar{\D} u_k (t,\ove{x})) : \bar{\D} v_k (\ove{x}).
			\end{align*}
		 	We argue as in Lemma~\ref{lem:convergence}.  By Assumption~\ref{assumptions} ~(ii)     and (iii) the energy density $W$ is quadratic near $Z$, i.e., there exists $\eta > 0$ such that  $ \vert DW  (Z + F) \vert  \le c |F|$ for all $F \in \R^{d \times 2^d}$ with $|F|\le \eta$. 	 On the other hand, given $F \in \R^{d \times 2^d}$ with  $|F| > \eta$, by Assumption~\ref{assumptions}~(v) we can bound $\vert DW  (Z + F) \vert \leq c |F|^2+ c \le c (1+ 1/\eta^2) |F|^2$. 			Combining  both estimates,  we conclude the existence of a constant $c> 0$ such that  
		$$\frac{1}{\delta_k}DW(Z + \delta_k \bar{\D} u_k(t,\cdot)) \leq \frac{1}{\delta_k} \big( c  \delta_k	\vert \bar{\D} u_k (t, \cdot)\vert + c \delta^2_k 	\vert \bar{\D} u_k (t, \cdot)\vert^2 \big) \le c(1 + \vert \bar{\D} u_k (t, \cdot)\vert^2)$$
		 a.e.\ for almost every $t\in (0,T)$. Here, the last step follows from Young's inequality. 			Since the discrete gradient is defined everywhere in $\ti{\Omega}$, we can replace the sum by the integral as in \eqref{eq:sum and int} and hence obtain   
			\begin{align*}
				\epsi_k^d \det A \sum_{\ove{x} \in  \L'_{\epsi} (\ti{\Omega})} \frac{1}{\delta_k} DW(Z + \delta_k \bar{\D} u_k (t, \ove{x})) : \bar{\D} v_k (\ove{x}) \leq c  \big(1+ \Vert \bar{\D} u_k (t, \cdot)\Vert^2_{ L^2(\ti{\Omega})} \big)  \Vert \bar{\D} v_k \Vert_{L^\infty(\ti{\Omega})}.
			\end{align*}
			Now, using \eqref{eq:L infty bound of disc grad}, $H^{\ell}_0(\Omega;\R^d) \subset W^{1,\infty}(\Omega;\R^d)$, $\sup_k \Vert \bar{\D} u_k  (t,\cdot)  \Vert_{L^2(\ti{\Omega})}  \le c$ (uniformly in $t$)  by \eqref{eq:uniform energy} and the  compactness from Theorem~\ref{thm:gamma convergence and compactness}, we arrive at
			\begin{equation*}
					\la \ove{\partial I_k(u_k(t))}, v \ra_{ H^{\ell}_0(\Omega)  }  \leq c \Vert \D  v \Vert_{L^\infty(\Omega)} \leq c \Vert v \Vert_{ H^{\ell}_0(\Omega)} \quad \text{for almost every $t \in (0,T)$.}
			\end{equation*}
			 In particular, by integrating in time we obtain (i).

  We now address (ii). Due to  
$$
\int_0^T   \left\Vert \rho \ove{\ddot{u}_k(t)} \right\Vert_{H^{-\ell}(\Omega)}^2 	\dd t \leq 2 	\int_0^T   \left\Vert \rho\ove{\ddot{u}_k(t)} + \ove{\partial I_k(u_k(t)) } \right\Vert_{H^{-\ell}(\Omega)}^2  +  2	\int_0^T   \left\Vert \ove{\partial I_k (u_k(t) ) } \right\Vert_{H^{-\ell}(\Omega)}^2
$$
		  and (i), it suffices to bound the first term on the right-hand side. To  this aim, for fixed  $t \in (0,T)$,   by   definitions \eqref{eq:ove ddot u}--\eqref{eq:lip op} and the interpretation of $\partial I_k(u_k)$ as scaled lattice  strain  \eqref{eq:defin disc grad as lattice disp}, we get  
		\begin{align*}
			 \left|	 \langle \rho\ove{\ddot{u}_k(t)} + \ove{\partial I_k(u_k(t)) }, v \ra_{H^{\ell}_0(\Omega)}  \right|^2  &=  \left| (\rho\ddot{u}_k(t) + \partial I_k(u_k(t)) , v_k )_k  \right|^2  \\
			 &\leq \Vert v_k \Vert_k^2  \, \Vert \rho\ddot{u}_k(t) + \partial I_k(u_k(t))  \Vert_k^2. 
		\end{align*}
 Since $\Vert v_k \Vert_k \leq \Vert v \Vert_{L^2(\Omega)}$ by Lemma \ref{lem:approx by lattice displacements} and $ H^{\ell}_0  (\Omega;\R^d) \subset L^2(\Omega;\R^d)$, we deduce
\begin{align*}
 \Vert \rho\ove{\ddot{u}_k(t)} + \ove{\partial I_k(u_k(t)) }\Vert_{H^{-\ell}(\Omega)}^2  \le c  \Vert \rho\ddot{u}_k(t) + \partial I_k(u_k(t))  \Vert_k^2
 \end{align*}
 for a.e.\ $t \in (0,T)$.
 Then, by \eqref{eq:unform bound on force}  we  conclude. 
				\end{proof}
		\sloppy
		 The lemma above implies the existence of $\xi \in L^2(0,T;H^{-\ell}(\Omega;\R^d))$ such that $ \ove{\partial I_k (u_k)} \wto\xi$ in $L^2(0,T;H^{-\ell}(\Omega;\R^d))$.  Lemma~\ref{lem:convergence} and strong convergence pointwise in time from Lemma~\ref{lem:Aubin Lions} allows us to identify $ \xi = \partial I(w)$.  Moreover,  there  exists  $\lambda \in L^2(0,T;H^{-\ell}(\Omega;\R^d))$ such that 
		\begin{align}\label{eq:lambda}
			\rho \ove{\ddot{u}_k} \wto \lambda \quad \text{in $ L^2(0,T;H^{-\ell}(\Omega;\R^d))$}.
		\end{align}
	 It remains to  identify    $\lambda = \rho \partial_{tt} w$. 	 To this end, let   $\varphi \in C_\comp^\infty ([0,T] \times \Omega; \R^d)$ be of the form   $\varphi(t,x) = \psi(t) v(x)$ for  $\psi \in C_\comp^\infty([0,T])$ and $v \in C_\comp^\infty (\Omega;\R^d)$. By \eqref{eq:lambda} we obtain 
		\begin{align}\label{eq:ddot uk 1}
			\int_0^T \la   \rho   \ove{\ddot{u}_k(t)}, \varphi(t) \ra_{H^{\ell}_0(\Omega)} \dd t \to 	\int_0^T \la \lambda(t), \varphi (t) \ra_{H^{\ell}_0(\Omega)} \dd t. 
				\end{align}
						On the other hand, writing $\varphi_k(t) = \psi(t) v_k$  for each $t \in [0,T]$, where $v_k$ are the scaled lattice displacements from Lemma~\ref{lem:approx by lattice displacements},  we have by \eqref{eq:ove ddot u} that
		\begin{align}\label{eq:ddot uk 2}
				&\int_0^T \la   \rho   \ove{\ddot{u}_k(t)}, \varphi (t) \ra_{H^{\ell}_0(\Omega)} \dd t = \int_0^T \big(  \rho    \ddot{u}_k(t), \varphi_k (t) \big)_k \dd t  =  \int_0^T  \ddot{\psi}(t) ( \rho   u_k(t),  v_k )_k  \dd t \nonumber  \\
				&\to \int_0^T  \ddot{\psi}(t)  (  \rho   w(t), v) {  \, \dd t  } =  \int_0^T    (  \rho   w(t), \ddot{\varphi}(t)) {  \, \dd t,  } 
		\end{align}
where the  convergence follows from Lemma~\ref{lem:conv of norms}  and Lemma \ref{lem:Aubin Lions}.     By the density of  linear combinations of  test functions of the above form  and \eqref{eq:ddot uk 1}--\eqref{eq:ddot uk 2} we get $ \lambda =   \rho  \partial_{tt} w$ in $L^2(0,T;H^{-\ell}(\Omega;\R^d))$. 
		\fussy
 From \eqref{eq:ddot uk 1} we can also deduce pointwise convergence of 	$\partial_t \ti{u}_k(t)$ for all $t \in [0,T]$. We formulate and prove this property here for later purposes.

	\begin{lemma}[Convergence of $\partial_t \ti{u}_k$]\label{lem:weak cont}
		Let $(u_k)_k$ be a sequence of solutions to the atomistic problem \eqref{eq:atomistic the second} and  let $\ti{w}$  be its limit from Lemma~\ref{lem:Aubin Lions}. Then
		\begin{equation*}
				\partial_t \ti{u}_k(t) \wto \partial_t  \ti{w}(t)  \text{  weakly in $L^2(\ti{\Omega};\R^d)$ for all $t \in [0,T]$}.
		\end{equation*}
	\end{lemma}	\begin{proof} 
	 	Similarly to \eqref{eq:ddot uk 1}--\eqref{eq:ddot uk 2},  we consider a test function $\varphi \in C^\infty ([0,T] \times \Omega; \R^d)$  of the form   $\varphi(t,x) = \psi(t) v(x)$, now  for  $v \in C_\comp^\infty (\Omega;\R^d)$ and  $\psi \in C^\infty([0,T])$ with $\psi(T) = 0$. As $\partial_{tt} w  \in L^2(0,T;H^{-\ell}(\Omega;\R^d))$, we get $\partial_{t} w  \in C([0,T];H^{-\ell}(\Omega;\R^d))$ and, for fixed $t \in [0,T]$, we obtain by integration by parts 
			\begin{align}\label{LLL1}
				&\int_t^T  \la  \partial_{tt} w(s), \varphi (s) \ra_{H^{\ell}_0(\Omega)} \dd s  = -  \int_t^T \big(      \partial_t w(s), \partial_t \varphi(s) \big) \dd s   - (   \partial_t w(t),  \varphi(t) ) 
		\end{align}
		where we also used   $\partial_t \ti{ w } \in L^2(0,T;L^2(\ti{\Omega};\R^d))$, see \eqref{eq:weak conv dot u}. In a similar fashion,  we get
				\begin{align}\label{LLL2}
				&\int_t^T \la     \ove{\ddot{u}_k(s)}, \varphi (s) \ra_{H^{\ell}_0(\Omega)} \dd s  =  -  \int_t^T  \big(   \dot{u}_k(s), (\partial_t\varphi)_k(s) \big)_k  \dd s  - (   \dot{u}_k(t),  \psi(t) v_k )_k,   
		\end{align}	
where	$(\partial_t\varphi)_k(s) = \partial_t\psi(s) v_k$  for each $s \in [0,T]$. By using \eqref{eq:apriori} and \eqref{eq:norms are the same}, we find a (not relabeled) subsequence and  $z(t) \in L^2(\ti{\Omega};\R^d)$ such that $\partial_t \ti{u}_k(t) \rightharpoonup  z(t)$ in $L^2(\ti{\Omega};\R^d)$. Then, by   \eqref{eq:weak conv dot u}, \eqref{eq:ddot uk 1} (with $\lambda =  \rho  \partial_{tt} w$), and  Lemma \ref{lem:conv of norms}, we pass to the limit in \eqref{LLL2}, 		 and we compare with \eqref{LLL1} to find $(   \partial_t w(t),  \varphi(t) )    = (   z(t),  \varphi(t) )$. By the arbitrariness of the test function, we get $z(t) = \partial_t \ti{w}(t)$. Uniqueness of the limit eventually shows that the whole sequence $\partial_t \ti{u}_k(t)$ converges to $\partial_t \ti{w}(t)$. This concludes the proof.  
\end{proof}

\subsubsection{Joint limit}\label{sec: XXX3}

		 From definitions \eqref{eq:ove ddot u}--\eqref{eq:lip op} and the interpretation of $\partial I_k(u_k)$ as scaled lattice  strain,    see  \eqref{eq:defin disc grad as lattice disp}, for each $v \in C^\infty_c(\Omega)$ we get $\la \rho\ove{\ddot{u}_k(t)} + \ove{\partial I_k(u_k(t)) }, v \ra_{H^{\ell}_0(\Omega)} = (\rho\ddot{u}_k(t) + \partial I_k(u_k(t)) , v_k )_k$. Moreover,  for test functions  $\varphi \in C_\comp^\infty ([0,T] \times \Omega; \R^d)$ of the form   $\varphi(t,x) = \psi(t) v(x)$,  with   $v \in C_\comp^\infty (\Omega;\R^d)$ and  $\psi \in C^\infty_\comp([0,T])$,  we have 
		\begin{align*}
		\int_0^T \big| ((\rho \ddot{u}_k &  + \partial I_k(u_k))^\sim(t),  \psi(t) v  ) - (\rho\ddot{u}_k(t) + \partial I_k(u_k(t)) ,  \psi(t)   v_k )_k \big|^2\, {\rm d} t \notag \\ & \le \Vert \ti{v}_k- v \Vert^2_{L^2(\ti{\Omega})}  \Vert \psi \Vert_{L^\infty(0,T)}  \int_0^T \Vert (\rho \ddot{u}_k + \partial I_k(u_k))^\sim(t) \Vert^2_{L^2(\ti{\Omega})} \,  {\rm d}t   \to 0 
		\end{align*}
		 by \eqref{eq:norms are the same}, \eqref{eq:unform bound on force},  Lemmas \ref{lem:approx by lattice displacements}--\ref{lem:conv of norms}, and H\"older's inequality. Hence, 	 recalling \eqref{eq:conv to alpha}  and the convergences established in Sections \ref{sec: XXX1}--\ref{sec: XXX2},    by uniqueness of weak limits  we conclude  that 		
		\begin{align}\label{identitt}
			\alpha = \rho \partial_{tt}  \ti{w}  + \partial I(\ti{w})  \quad  \text{in $L^2(0,T;L^2(   \ti{\Omega};    \R^d))$.}
		\end{align}

		\subsection{Proof of   the  first part of Theorem~\ref{thm:main}:  convergence of solutions}\label{subsec:conclusion thm 1}
	
			Collecting the results from Lemma~\ref{lem:Aubin Lions},  \eqref{eq:weak conv dot u}, \eqref{eq:conv to alpha},  and \eqref{identitt}  we conclude   the existence of a (non relabeled) subsequence of  solutions $(u_k)_k$ to the atomistic equation of motion \eqref{eq:atomistic the second}, satisfying
				\begin{align}
					& \ti{u}_k(t) \to \ti{w}(t) \quad &&\text{in $L^2(\ti{\Omega};\R^d)$  for all $t \in [0,T]$,}\label{auchnummer} \\
					& \nu \partial_{t} \ti{u}_k \wto \nu \partial_{t}  \ti{w}  \quad &&\text{in  } L^2(0,T;L^2(\ti{\Omega};\R^d)), \label{eq:weak conv to sol 1} \\
					&\rho \partial_{tt} \ti{u}_k + (\partial I_k(u_k))^\sim \wto \rho \partial_{tt}  \ti{w}  +  \partial I (\ti{w})  \quad &&\text{in  } L^2(0,T;L^2(\ti{\Omega};\R^d)).    \label{eq:weak conv to sol 2}
				\end{align}
			In particular,  testing \eqref{eq:atomistic the second} with $\psi(t) v_k$, where $\psi \in C^\infty([0,T])$ and   $v_k$ is the associated scaled lattice displacements  for some $v \in C^\infty_\comp(\Omega)$, and passing to the limit $k \to \infty$, we find that 		  $w$ is  a weak solution of the continuum momentum equation.  Moreover, the  initial conditions are satisfied.  In fact,  by Lemma~\ref{lem:Aubin Lions} we have  $u_k(0) \to w(0)$ in $L^2(\ti{\Omega};\R^d)$,  and Lemma~\ref{lem:weak cont} implies   $\partial_t \ti{u}_k (0) \wto \partial_t w(0)$ in $L^2(\ti{\Omega};\R^d)$. This along with  $ \ti{u}_k (0, \cdot ) \to w^0$  and  $\partial_t \ti{u}_k (0, \cdot ) \to w^1$ in $L^2(\ti{\Omega};\R^d)$ by assumption   yields the claim.  By $w\in  L^\infty(0,T;H^1_0(\Omega;\R^d))$ (see Lemma \ref{lem:Aubin Lions}), \eqref{eq: forrrlater}, and    comparison, we have $\partial_{tt} w \in L^2(0,T; H^{-1}(\Omega;\R^d))$.

			 As the continuum problem \eqref{eq:continuous problem}  admits a unique weak solution $w$,  the entire sequence $(u_k)_k$ of solutions to the atomistic problem AC-converges.  \qed

				The convergence  given in \eqref{eq:weak conv to sol 1} and \eqref{eq:weak conv to sol 2} can be improved under the stronger assumption on the initial data in the second part of Theorem~\ref{thm:main}, see Corollary~\ref{cor: cor} below for details.

	\section{Energy convergence}\label{sec:proof 2}

		We devote this section to the proof  the energy-convergence statement  in  Theorem~\ref{thm:main}. In Section~\ref{subsec:lsc}, we prove inequalities   \eqref{eq:liminf for time derivative}--\eqref{eq:liminf for grad of energy}.  The proof of Theorem~\ref{thm:main} is then concluded in Section~\ref{subsec:concluding the proof}.

	\subsection{lim\,inf-inequalities}\label{subsec:lsc}
	 Having already established the convergences \eqref{auchnummer}--\eqref{eq:weak conv to sol 2},  we now simply use the weak lower semicontinuity of norms to establish the inequalities.

		   Inequality \eqref{eq:liminf for time derivative} follows from Lemma~\ref{lem:weak cont},  \eqref{eq:norms are the same},  and the weak lower semicontinuity of norms.   Then, \eqref{eq:liminf for integral of time derivative} follows by integration and 	 Fatou's Lemma.

		  For \eqref{eq:zoom12}, we use \eqref{auchnummer}, Lemma \ref{lem:conv of norms}, and the $\Gamma$-lim inf inequality of Theorem \ref{thm:gamma convergence and compactness}. 
 	 
 	  	 Eventually,  inequality~(\ref{eq:liminf for grad of energy}) follows   from \eqref{eq:conv to alpha},  \eqref{identitt},  and weak lower semicontinuity of norms,  again recalling relation \eqref{eq:norms are the same}.

  		\begin{remark}[Purely viscous case $\rho =0$] \label{rem:inertia grad in L2}
  		 Let us address the analog of the inequalities \eqref{eq:liminf for time derivative}--\eqref{eq:liminf for grad of energy} in the purely  viscous setting.  First, \eqref{eq:liminf for time derivative} is not needed and  \eqref{eq:zoom12} remains unchanged. Inequality \eqref{eq:liminf for integral of time derivative} still holds, but for its derivation we do not resort to \eqref{eq:liminf for time derivative}, but simply use  \eqref{eq:weak conv dot u}, lower semicontinuity of norms, and   \eqref{eq:norms are the same}.   		
  		
 The only inequality which essentially changes is  \eqref{eq:liminf for grad of energy}: it is replaced    by  
  		\begin{align}\label{rhogleich0} 
  			\Vert \div (\C: \D^s w(t, \cdot) ) \Vert_{L^2(\Omega)} \leq \liminf_{ k \to \infty} \Vert \partial I_k(u_k(t, \cdot)) \Vert_k  \quad \text{for a.e.\ $t \in (0,T)$}. 
  		\end{align} 
  		 Once \eqref{rhogleich0} is shown, we can use the uniform bound on $\int_0^T \Vert \partial I_k (u_k(s)) \Vert_k^2 \dd s$, see \eqref{eq:EDIE} (for $\rho = 0$) and Fatou's Lemma to conclude  that  $\partial I (w) \in L^2(0,T;L^2(\Omega;\R^d))$,   cf.\ \eqref{eq:derivative of limiting energy}.    		  
  		
  		Let us show \eqref{rhogleich0}. By summation by parts, see Lemma~\ref{lem:summation by parts},  $\la \partial I_k (u_k(t)) , v_k \ra$ can be expressed via the atomistic inner product $( \cdot, \cdot)_k$ as 
  		\begin{align*}
  			\epsi_k^d \det A \sum_{\ove{x} \in  \L'_{\epsi_k} (\ti{\Omega})} \frac{1}{\delta_k} DW(Z + \delta_k \bar{\D} u_k (t, \ove{x})) \cdot \bar{\D} v_k(\ove{x}) = (f_k(t), v_k)_k 
  		\end{align*}
  		for any $v \in C_c^\infty(\Omega;\R^d)$ and the corresponding $v_k$ from Lemma~{\ref{lem:approx by lattice displacements}},  	where the specific form of $f_k(t)$ can be found in \eqref{eq:conjugate operator} below.  It is not restrictive to assume that $\liminf_{k \to \infty} \Vert f_k (t) \Vert_k < \infty$ as otherwise  \eqref{rhogleich0} 	 is trivial.   Then, by   Lemma \ref{lem:approx by lattice displacements}   we get
  		$$\liminf_{k \to \infty}(f_k(t), v_k)_k \le \liminf_{k \to \infty} \Vert f_k(t) \Vert_k \Vert v_k \Vert_k \le c \Vert v \Vert_{L^2(\Omega)}.  $$
  		As 	  Lemma~\ref{lem:convergence} implies 
  		$$ 	 \langle   \partial I( w(t) ) , v 	 \rangle  =  \int_\Omega   (\D  w (t) \cdot Z) :\mathbb{H}: (\D v \cdot Z) \dd x  = \lim_{k\to \infty} (f_k(t), v_k)_k,  $$
  		we get that $\partial I( w(t) )$ is bounded in $L^2(\Omega;\R^d)$. Then,  \eqref{eq:derivative of limiting energy} yields 
  		$$ (f_k(t), v_k)_k \to ( \partial I(w(t)) , v ) = 	 - \int_\Omega  \div (\C : \D^s  w (t))\cdot v \dd x = (f(t),v), $$		
  		where $f(t) := - \div (\C : \D^s  w (t) )$.  Denoting   	by  $\ti{f}_k(t)$    the piecewise constant   interpolation of $f_k(t)$  introduced in \eqref{eq:pw const int},  we  claim   $\ti{f}_k(t) \rightharpoonup f(t)$ weakly in $L^2(\Omega;\R^d)$. This  would then show   $	\Vert f(t) \Vert_{L^2(\Omega)} \leq 	 \liminf_{k \to \infty}  \Vert \ti{f}_k(t) \Vert_{L^2( \ti\Omega)} $ which, along with \eqref{eq:norms are the same}, concludes the proof of \eqref{rhogleich0}. 
  		
  		 We are left with proving that  $\ti{f}_k(t) \rightharpoonup f(t)$ weakly in $L^2(\ti{\Omega};\R^d)$.  Given $v \in C^\infty_\comp(\Omega)$, we note   that $| (\ti{f}_k(t)  , v) - (\partial I(w(t)), v)| \leq | (\ti{f}_k(t)  , v) - (f_k(t), v_k)_k | + | (f_k(t), v_k)_k- (\partial I(w(t)) , v) | $ and hence by \eqref{eq:sum and int} it is enough to check that $| (\ti{f}_k(t)  , v) - (\ti{f}_k(t), \ti{v}_k) |$ vanishes as $k \to \infty$. Indeed, by H\"{o}lder's inequality  and \eqref{eq:norms are the same}  we have
  		\begin{align*}
  			| (\ti{f}_k(t)  , v) - (\ti{f}_k(t), \ti{v}_k) |^2 \leq \Vert v - \ti{v}_k \Vert_{L^2(\ti{\Omega})}^2 \Vert f_k \Vert_k^2 =  \Vert v - \ti{v}_k \Vert_{L^2(\ti{\Omega})}^2 \Vert \partial I_k (u_k(s)) \Vert_k^2  \to 0,  
  		\end{align*}
 where the convergence follows from $\liminf_{k \to \infty} \Vert f_k (t) \Vert_k < \infty$ and Lemma \ref{lem:conv of norms}. 
  	\end{remark}

		\subsection{Conclusion of the proof of Theorem~\ref{thm:main}}\label{subsec:concluding the proof}
		
		In order to follow the evolutionary $\Gamma$\nobreakdash-convergence approach, the chain rule for $\rho \partial_{tt} w + \partial I (w)$ is needed. For  our  weak solutions, however,  with regularity  $\partial_{tt} w \in L^2(0,T;H^{-1}(\Omega; \R^d))$ and $\partial_{t} w \in L^2(0,T;L^2(\Omega;\R^d))$,  the  required  chain rule in $L^2$ is not available.  		 By assuming  additional regularity on the initial data, we can prove  the existence of  strong solutions to the continuum momentum equation \eqref{eq:continuous problem},  which are sufficiently regular  such that  a chain rule in $L^2$ holds. 
		
			 Indeed, standard arguments, see, e.g., \cite[Prop.~6.3.2 on p.~204]{kruzikMathematicalMethodsContinuum2019}, give the following. 
		
		\begin{lem}[Strong solutions to \eqref{eq:continuous problem} and the chain rule]\label{lem:chain rule}
			Let $w^0 \in H^2(\Omega;\R^d)$ and $w^1 \in H^1(\Omega;\R^d)$. Then,  the unique weak solution $w$  to the continuum momentum equation \eqref{eq:continuous problem}  is  a strong solution  with
			\begin{align*}
				&w \in L^{\infty}(0,T;  H^2(\Omega;\R^d) ), \\
				&\partial_{t} w \in L^{\infty}(0,T; H^1_0(\Omega;\R^d)), \\
				&\partial_{tt} w  \in L^{\infty}(0,T; L^2(\Omega;\R^d)).
			\end{align*}
			In particular, $\partial I(w) \in L^2(0,T;L^2(\Omega;\R^d))$  by comparison  and the chain rule
			\begin{align}\label{eq:chain rule}
		\frac{\dd}{\dd s} \left( \frac{\rho}{2} \Vert \partial_t w (s, \cdot) \Vert_{L^2(\Omega)}^2 + I(w(s)) \right) =	 \big( \rho \partial_{tt}  w   (s, \cdot ) +  \partial I(w   (s, \cdot )), \partial_t w (s, \cdot) \big)
			\end{align} 
		holds	for almost every $s \in (0,T)$.
		\end{lem}

\begin{proof}[Proof of Theorem~\ref{thm:main}]   To obtain energy convergence, our idea is to  pass to the limit in the atomistic energy-dissipation~inertia equalities by arguing along the   lines of {\sc Sandier \& Serfaty} \cite{ss}.  			As   $u_k$ is a solution to the  atomistic equation of motion  \eqref{eq:atomistic the second},
			  \eqref{eq:EDIE} holds, i.e.,  
			\begin{align}\label{eq: again and again}
				I_{k} &(u_k(0)) - I_k (u_k(t)) +  \frac\rho2 \Vert \dot{u}_k (0)\Vert_k^2 - \frac\rho2 \Vert \dot{u}_k (t)\Vert_k^2 \notag \\
							&= \frac\nu2 \int_0^t  \Vert \dot{u}_k (s,\cdot ) \Vert_k^2 \dd s + \frac{1}{2\nu} \int_0^t  \Vert \rho \ddot{u}_k(s, \cdot ) +  \partial I_k (u_k(s, \cdot )) \Vert_k^2 \dd s\quad  \forall t \in [0,T]. 
			\end{align}
			 By   Lemma~\ref{lem:Aubin Lions}   we find $ w(t, \cdot )  \in  H^1_0  (\Omega; \R^d)$  such that $u_{k}(t, \cdot) \acto  w  (t, \cdot)$  for  all  $t \in [0,T]$.   	Therefore,  by means of inequalities  \eqref{eq:liminf for time derivative}--\eqref{eq:liminf for grad of energy},  $I_{k}( u_{k} ) \to I( w^0 )$,  and  $ \dot{u}_{k}(0)  \acto  \partial_t w(0) $,  we get 
		\begin{align}\label{eq: the 1}
			&I( w (0)) - I( w  (t)) +  \frac\rho2 \Vert  \partial_t w (0, \cdot) \Vert_{L^2(\Omega)}^2 - \frac\rho2 \Vert \partial_t w(t, \cdot) \Vert_{L^2(\Omega)}^2 \notag \\
			&\geq\liminf_{k \to \infty} \big( I_k(u_k(0))-I_k (u_k (t) )  \big) +  \liminf_{ k \to \infty} \big( \frac\rho2 \Vert \dot{u}_k (0)\Vert_k^2 - \frac\rho2 \Vert \dot{u}_k (t)\Vert_k^2 \big)\notag\\
			&\geq \liminf_{k \to \infty}   \frac\nu2  \int_0^t \Vert \dot{u}_k (s,\cdot ) \Vert_k^2 \dd s +  \liminf_{k \to \infty} \frac{1}{2\nu} \int_0^t  \Vert  \rho \ddot{u}_k (s, \cdot) +  \partial I_k (u_k (s, \cdot)) \Vert_k^2 \dd s \nonumber \\
			&\geq   \frac\nu2  \int_0^t \Vert\partial_t  w   (s, \cdot ) \Vert_{L^2(\Omega)}^2 \dd s + \frac{1}{2\nu} \int_0^t \Vert \rho \partial_{tt} w (s, \cdot) +  \partial I ( w   (s, \cdot ) ) \Vert_{L^2(\Omega)}^2 \dd s.
		\end{align}	
		 Using Young's inequality we can further estimate 
		\begin{align}\label{eq: the 2}
			 &\frac\nu2  \int_0^t \Vert\partial_t  w   (s, \cdot ) \Vert_{L^2(\Omega)}^2 \dd s + \frac{1}{2\nu} \int_0^t \Vert \rho \partial_{tt} w (s, \cdot) +  \partial I ( w   (s, \cdot ) ) \Vert_{L^2(\Omega)}^2 \dd s \nonumber \\
			&\geq - \int_0^t (  \rho \partial_{tt}  w   (s, \cdot ) +  \partial I(w   (s, \cdot )), \partial_t w (s, \cdot) ) \dd s. 
		\end{align}
		Hence, from the chain rule \eqref{eq:chain rule},  the right-hand side of \eqref{eq: the 2} equals  
		\begin{align}\label{eq: the 3}
			-&\int_0^t \frac{\dd}{\dd s} \left( \frac{\rho}{2} \Vert \partial_t w (s, \cdot) \Vert_{L^2(\Omega)}^2 + I(w(s)) \right) \dd s\notag\\
			&= I( w  (0)) - I( w  (t)) +  \frac\rho2 \Vert  \partial_t w (0, \cdot) \Vert_{L^2(\Omega)}^2 - \frac\rho2 \Vert \partial_t w(t, \cdot) \Vert_{L^2(\Omega)}^2. 
		\end{align}
		 Collecting \eqref{eq: the 1}--\eqref{eq: the 3} we thus conclude that  all inequalities  above are actually equalities.  This is exactly the  energy-dissipation-inertia equality  corresponding to the continuum equation of motion  \eqref{eq:continuous problem}, i.e., we have  rederived that  $w$ is a solution to \eqref{eq:continuous problem}. In contrast to the reasoning in Section \ref{subsec:conclusion thm 1}, the present approach has the advantage of  yielding energy convergence   $I_k(u_k(t)) \to I( w  (t))$ for all $t \in (0,T)$.	Eventually,  energy convergence  entails the  strong convergence of  gradients, namely $\bar{\D} u_k(t) \to \D w(t)\cdot Z$ strongly in $L^2(\Omega;\R^d)$ for  all  $t \in (0,T)$,  see Remark~\ref{rem:bdry and Cauchy Born}. 
\end{proof}

 By  inspecting the above  argument,   the convergence  can be strengthened as follows.  
  
		\begin{corollary}[Enhanced convergence]\label{cor: cor}
			Let $u_k$ and $ w $ solve  \eqref{eq:atomistic the second}  and \eqref{eq:continuous problem} respectively. Under the assumptions of  the second part of   Theorem~\emph{\ref{thm:main}}, we have
			\begin{align*}
				& \mathrm{(i)}	&&\int_0^T	\Vert   \partial_t \ti{u}_k  (t, \cdot ) - \partial_t \ti{w} (t, \cdot ) \Vert^2_{ L^2(\ti{\Omega})} \dd t\to 0, \\
				&   \mathrm{(ii)}	&& \int_0^T	\Vert   (\rho \ddot{u}_k(t) + \partial I_k(u_k(t)))^\sim     -  ( \rho \partial_{tt} \ti{w}(t, \cdot ) +    \partial I ( \ti{w} ))(t,\cdot ) \Vert^2_{ L^2(\ti{\Omega})}  \dd t\to 0.   \\
			\end{align*}
		\end{corollary}
	\begin{proof}
 Observe that weak convergence has already been established in \eqref{eq:weak conv to sol 1}--\eqref{eq:weak conv to sol 2}. Strong convergence in $L^2(0,T;L^2(\ti{\Omega};\R^d))$ then follows by checking that also the norms convergence. The latter is a consequence  of \eqref{eq:norms are the same} and  the evolutive $\Gamma$\nobreakdash-convergence in proof of Theorem~\ref{thm:main}, again observing that all inequalities are actually equalities.  
	\end{proof}

	\begin{remark}[Purely viscous case $\rho = 0$]\label{rem:purely viscous setting strong solutions}
 In the case $\rho = 0$, we first note that we still can identify a limit $w$, see Remark \ref{rem:inertia grad in L2-again}.  Remark~\ref{rem:inertia grad in L2} yields  $\div(\C :	\D^s w) \in L^2(0,T;L^2(\Omega;\R^d))$. This ensures the validity of the $L^2$ chain rule under no additional regularity assumption on the initial data.    Therefore, existence of solutions needs not be proven beforehand, but follows from passing to the limit in the  energy-dissipation  equation \eqref{eq: again and again} for $\rho = 0$,  by resorting to the inequalities  \eqref{eq:liminf for integral of time derivative},  \eqref{eq:zoom12},     \eqref{rhogleich0},  and using Fatou's Lemma.  
	Eventually, standard elliptic regularity  also  implies 	 $w \in	L^2(0,T;H^2(\Omega;\R^d))$,  i.e., $w$ is the unique   strong solutions of the continuum  momentum equation. 
\end{remark}

	\section*{Acknowledgments} 
This work was supported by the DFG project FR 4083/3-1, by the
Deutsche Forschungsgemeinschaft (DFG, German Research Foundation)
under Germany's Excellence Strategy EXC 2044 -390685587, Mathematics
M\"unster: Dynamics--Geometry--Structure,   by the Austrian
Science Fund (FWF) projects F\,65, I\,4354, I\,5149, and P\,32788.   The authors would like to thank \name{Marco Cicalese}, \name{Mark Peletier}, and \name{Florian Theil}  for some interesting comments on a previous version of the paper.

		\appendix
		\section{} \label{sec:appendix}

		\begin{lemma}[Summation by parts]\label{lem:summation by parts}
	Let $  g  = (g_1, \dots, g_{2^d}) \colon \L'_{\epsi}(\Omega) \to \R^{d \times 2^d}$, i.e., $g_i(\ove{x}) \in \R^d$ for $i = 1, \dots, 2^d$,  $\bar{x} \in \L'_{\epsi}(\Omega)$,  and let $v \in \A_\epsi$ for $\epsi > 0$ small enough.  Then, the following \df{summation by parts formula} holds
			\begin{equation}\label{eq:summation by parts}
				\sum_{\ove{x} \in  \L'_\epsi(\ti{\Omega})}   g(\ove{x}) \colon   \bar{\D} v (\ove{x}) = -\sum_{x \in \L_\epsi \cap  {\Omega}} \bar{\D}^* g(x) \cdot v(x),
			\end{equation}
			where $\bar{\D} v $ is the discrete gradient defined in \eqref{eq:def of disc grad} and  $\bar{\D}^*$ denotes the conjugate operator to $\bar{\D}$ given by 
			\begin{equation}\label{eq:conjugate operator}
				\bar{\D}^* g (x) =  \frac{1}{\epsi}  \sum_{i = 1}^{2^d} \left(  - g_i (x - \epsi z_i)+ \frac{1}{2^d} \sum_{j =1}^{2^d} g_i (x - \epsi z_j) \right)  \quad   \forall \, x \in \L_\epsi \cap {\Omega}.
			\end{equation}
		\end{lemma}
	\begin{proof}
		Recall the definition of $\bar{\D} v (\ove{x}) = \frac{1}{\epsi} \bigl(v (\ove{x} + \epsi z_1)- \bar{v}   (\ove{x}),  \dots , v (\ove{x} + \epsi z_{2^d})- \bar{v}   (\ove{x})   \bigr)$, where $\bar{v} (\ove{x})= \frac{1}{2^d} \sum_{j = 1}^{2^d} v(\ove{x} + \epsi z_j)$. We abbreviate $v_i (\ove{x}) := v(\ove{x} + \epsi z_i)$ in the following. 
		We split the left hand side of \eqref{eq:summation by parts} as
		\begin{align*}
			\sum_{\ove{x} \in \L'_\epsi(\ti{\Omega})} g(\ove{x}) \colon \bar{\D} v(\ove{x}) = \underbrace{ \frac{1}{\epsi}\sum_{\ove{x} \in \L'_\epsi ( \ti{\Omega} ) } g(\ove{x}) \colon (v_1  (\ove{x}),  \dots,  v_{2^d}  (\ove{x})  )   }_{:= I} - \underbrace{\frac{1}{\epsi} \sum_{\ove{x} \in \L'_\epsi ( \ti{\Omega} ) } g(\ove{x}) \colon  (\bar{v}  (\ove{x}),  \dots, \bar{v} (\ove{x}) )   }_{:= II}.
		\end{align*}
		For the first term we obtain
		\begin{align*}
			I & = \frac{1}{\epsi}\sum_{\ove{x} \in \L'_\epsi ( \ti{\Omega} ) } \sum_{i = 1}^{2^d} g_i (\ove{x}) \cdot v (\ove{x} + \epsi z_i)  =\frac{1}{\epsi} \sum_{x \in \L_\epsi \cap  {\Omega}  } v(x) \cdot \left( \sum_{i = 1}^{2^d} g_i (x - \epsi z_i) \right),
		\end{align*}
		which follows by the change of variables $\ove{x} + \epsi z_i = x$ and by noting that $v(x) = 0$ for $x \notin \Omega$. The second term can be  handled analogously as
		\begin{align*}
			II &= \frac{1}{\epsi}\sum_{\ove{x} \in \L'_\epsi ( \ti{\Omega} ) }  \sum_{i = 1}^{2^d} g_i (\ove{x}) \cdot \bar{v}(\ove{x}) = \frac{1}{\epsi}\sum_{\ove{x} \in \L'_\epsi (  \ti{\Omega} ) } \sum_{i= 1}^{2^d} g_i (\ove{x})  \cdot \left(    \frac{1}{2^d} \sum_{j = 1}^{2^d} v (\ove{x} + \epsi z_j)     \right) \\
			&= \frac{1}{\epsi}\sum_{x \in \L_\epsi \cap  {\Omega}  } v(x) \cdot \left( \sum_{i =1}^{2^d} \frac{1}{2^d} \sum_{j =1}^{2^d} g_i (x - \epsi z_j) \right),
		\end{align*}	
		with the change of variables as above. 
		\end{proof}

\end{document}